\documentclass[12pt]{article}
\usepackage{latexsym}
\usepackage{amssymb}
\usepackage{amsthm}
\usepackage{epsfig}
\usepackage{graphicx}
\usepackage{mathtools} 
\usepackage{enumitem}
\usepackage{xcolor}
\usepackage{easyReview}
\usepackage{hyperref}
\usepackage{algorithm}
\usepackage[noend]{algpseudocode}
\makeatletter
\def\BState{\State\hskip-\ALG@thistlm}
\makeatother

\DeclareMathOperator*{\argmin}{arg\,min}
\hypersetup{
    linkcolor=blue,  
    linkbordercolor=red,
    linktoc=all,     
    linktocpage
}
\advance \topmargin by -\headheight
\advance \topmargin by -\headsep
\evensidemargin \oddsidemargin
\allowdisplaybreaks
\begin{document}
\title{Optimal Interference Signal for Masking an Acoustic Source}
\author{\\ Hongyun Wang
\\
Department of Applied Mathematics \\
University of California, Santa Cruz, CA 95064, USA\\
\\
Hong Zhou \\ Department of Applied Mathematics\\
Naval Postgraduate School, Monterey, CA 93943, USA \\ \\
}

\maketitle
\thispagestyle{empty}

\pagestyle{plain}

\begin{abstract}
In an environment where acoustic privacy or deliberate signal obfuscation is desired, it is necessary to mask the acoustic signature generated in essential operations. We consider the problem of masking the effect of an acoustic source in a target region where possible detection sensors are located. 
Masking is achieved by placing interference signals near the acoustic source. 
We introduce a theoretical and computational framework for designing such interference signals with the goal of minimizing the residual amplitude in the target region. 
For the three-dimensional (3D) forced wave equation with spherical symmetry, 
we derive analytical quasi-steady periodic solutions for several canonical cases. 
We examine the phenomenon of self-masking where an acoustic source with certain spatial forcing profile masks itself from detection outside its forcing footprint. 
We then use superposition of spherically symmetric solutions to investigate masking in a given target region. 
We analyze and optimize the performance of using one or two point-forces 
deployed near the acoustic source for masking in the target region. 
For the general case where the spatial forcing profile of the acoustic source lacks spherical symmetry, 
we develop an efficient numerical method for solving the 3D wave equation. 
Potential applications of this work include undersea acoustic communication security, 
undersea vehicles stealth, and protection against acoustic surveillance.
\end{abstract}


\clearpage
\renewcommand*\contentsname{Table of contents}
\tableofcontents
\clearpage

\section{Introduction}
Controlling acoustic environments is an area of growing interest across disciplines such as architectural acoustics, audio security, and acoustic surveillance \cite{Kuttruff_2016, Yang_2017, Edelmann_2011, Ding_2024}. A particularly challenging problem in this domain is how to obscure or mask an acoustic source so that it becomes undetectable or unintelligible to a sensor. Traditional approaches rely on physical shielding or broadband noise injection, both of which have notable limitations related to effectiveness, energy efficiency, and intrusiveness \cite{Li_2016}.

Recent developments in signal processing and optimization have enabled more advanced and targeted acoustic masking techniques. In particular, the generation of interference signals that are carefully engineered to effectively cancel or perceptually mask an acoustic source has shown promise \cite{Hyde_2023, Bassett_2024}.
These methods have the potential to achieve truly effective masking with a lower energy cost and minimal disturbance to surrounding environments \cite{Jiang_2018, Liu_2019}.

In this paper, we address the problem of designing an optimal interference signal to mask the effect of a known acoustic source in a given region of sensors. 
We formulate this task as an optimization problem, with the objective of minimizing the perceptibility or detectability of the source signal by sensors in a specified region. 

The remainder of the paper is organized as follows. Section 2 introduces the mathematical formulation of the problem, including the analytical solution of the 3D wave equation under the assumption of spherical symmetry. Section 3 illustrates several examples of acoustic masking based on superposition
of spherically symmetric solutions. 
Section 4 presents an efficient numerical method for simulating the wave propagation in the 3D infinite space in the case where the spatial forcing profile of the acoustic source lacks spherical symmetry. 
Section 5 discusses the directions for future research and concludes the paper.

\section{Analytical Solution of the 3D Wave Equation with Spherical Symmetry}
\label{wav_sph_sym}
In this section, we derive solutions to the three-dimensional (3D) forced wave equation
under spherical symmetry. This approach is based on Duhamel's principle, which is equivalent to employing a
Green's function in the time dimension.
We begin by solving the special case in which the acoustic source is a temporal impulse forcing at an arbitrary time. The solution for a general forcing function is then constructed by superposing (or more precisely integrating) these impulse responses, yielding a time integral representation. 
The time integral solution contains the spatial forcing profile of the acoustic source. In general, this integral requires a numerical evaluation. In the special case where the spatial forcing profile is
a truncated $\text{sinc}(\cdot)$ function, this time integral solution has a closed-form expressions in terms of elementary trigonometric functions.

\subsection{The 3D Wave Equation with Spherical Symmetry}
We study acoustic wave propagation in the unbounded three-dimensional space $\mathbb{R}^3$. 
The governing equation is the wave equation with forcing: 
\begin{equation}
\begin{dcases}
u_{tt}(\vec{x}, t) = c^2 \nabla_{\vec{x}}^2 u(\vec{x}, t) 
+ F(\vec{x}, t) \\[1ex]
u(\vec{x}, 0) = u_0(\vec{x}), \quad u_t(\vec{x}, 0) = v_0(\vec{x})
\end{dcases}
\label{wav_eq1}
\end{equation}
where $u$ is the acoustic pressure at position ${\vec x}$ and time $t$; $c$ is the speed of acoustic wave in the medium under consideration; $\nabla^2$ is the Laplacian operator, and
$F$ represents the forcing of the acoustic source.
We consider the case where the force term and initial conditions satisfy the following properties: 
\begin{itemize}
\item The force term exhibits separable dependences on $\vec{x}$ and $t$: 
$F(\vec{x}, t) = \sigma(\vec{x}) f(t)$; 
\item 
The force term is spatially spherically symmetric:
$F(\vec{x}, t) = \sigma(r) f(t), \; r \equiv |\vec{x}|$ ; 
\item 
The initial conditions are zero: $u_0(\vec{x}) \equiv 0, \; v_0(\vec{x}) \equiv 0$. 
\end{itemize}
Under these assumptions, initial value problem (IVP) \eqref{wav_eq1} is spherically symmetric. 
We therefore express \eqref{wav_eq1} in spherical coordinates $(r, \theta, \phi)$.
Due to the spherical symmetry, the solution $u$ depends only on $(r, t)$, and is independent of 
$(\theta, \phi)$. Consequently, $u(r, t)$ satisfies the initial boundary value problem (IBVP) below. 
\begin{equation}
\begin{dcases}
u_{tt}(r, t) = c^2 \frac{1}{r^2}\Big(r^2 u_r(r, t) \Big)_r + \sigma(r) f(t), 
\quad r \in (0, \infty), \; t > 0\\[1ex]
u(0, t) \;\; \text{is finite} \\
u(r, 0) = 0, \quad u_t(r, 0) = 0
\end{dcases}
\label{wav_eq2}
\end{equation}

\subsection{Duhamel's Principle for Solving the Initial Boundary Value Problem \eqref{wav_eq2}}
Let $u^\text{(IV)}(r, t)$ denote the solution of \eqref{wav_eq2} in the special case where the forcing term is a temporal impulse at time 0: 
\begin{equation}
F(r, t)= \sigma(r) \delta(t-0)
\label{sp_case_0}
\end{equation}
The sole effect of this impulse force is to instantaneously increase the velocity from 
$u_t(r, 0^{-}) = 0$ to $u_t(r, 0^{+}) = \sigma(r)$. 
As a result, $u^\text{(IV)}(r, t)$ satisfies IBVP \eqref{wav_eq2} 
with zero forcing term and initial velocity updated to $u_t(r, 0^{+}) = \sigma(r)$, which leads 
to an IBVP on the homogeneous wave equation:
\begin{equation}
\begin{dcases}
u_{tt}(r, t) = c^2 \frac{1}{r^2} \Big(r^2 u_r(r, t) \Big)_r, 
\quad r \in (0, \infty), \; t > 0\\[1ex]
u(0, t) \;\; \text{is finite} \\
u(r, 0) = 0, \quad u_t(r, 0) = \sigma(r)
\end{dcases}
\label{wav_eq3}
\end{equation}
For simplicity and conciseness, we will use the generic notation $u$ to denote the unknown function 
in all IBVPs. For example, IBVP \eqref{wav_eq3} is written in terms of $u(r, t)$ although it is the governing IBVP for $u^\text{(IV)}(r, t)$. 
The solution of IBVP \eqref{wav_eq3} gives the value of $u^\text{(IV)}$ for $t \ge 0$. 
Before the start of forcing, the acoustic pressure is zero. Thus, we set $u^\text{(IV)}(r, t) = 0$
for $ t < 0$. After this extension, $u^\text{(IV)}(r, t)$ is defined for $t \in (-\infty, +\infty)$. 
We now switch to a special case that is slightly more general than \eqref{sp_case_0}.
\begin{equation}
F(r, t)= \sigma(r) \delta(t-t')
\label{sp_case_1}
\end{equation}
Then, the solution to IBVP \eqref{wav_eq2} with forcing term \eqref{sp_case_1} is
given by $u^\text{(IV)}(r, t-t')$, which is nonzero only for $t > t'$. 
In IBVP \eqref{wav_eq2}, the general forcing $\sigma(r) f(t)$ can be viewed as 
a superposition of impulse forcing of the form given in \eqref{sp_case_1}. 
Specifically, we write
\begin{equation}
 \sigma(r) f(t) = \int_0^{\infty} f(t') \Big(\sigma(r) \delta(t-t')\Big) dt' 
 \end{equation}
It follows that the solution to \eqref{wav_eq2} for the general forcing  
is a superposition of $u^\text{(IV)}(r, t-t')$. 
\begin{equation}
\boxed{\quad u(r, t) = \int_0^{\infty} f(t') u^\text{(IV)}(r, t-t') dt' 
= \int_0^t f(t-\tau) u^\text{(IV)}(r, \tau) d\tau \quad}
\label{Duhamel}
\end{equation}
To calculate $u(r, t)$ in \eqref{Duhamel}, we need to solve IBVP \eqref{wav_eq3}
to obtain $u^\text{(IV)}(r, t)$. We now employ a transformation to recast
\eqref{wav_eq3} as an IBVP on the one-dimensional (1D) wave equation. 
We introduce the transformation $w(r, t) = r u(r, t)$
and derive the governing IBVP for $w(r, t)$. 
First, we express the derivatives of 
$u(r, t)$ in terms of those of $w(r, t) $. 
\begin{equation}
\begin{dcases}
u(r, t) = \frac{1}{r}w(r, t), \quad 
u_t = \frac{1}{r} w_t, \qquad u_{tt} = \frac{1}{r} w_{tt} \\[1ex]
u_r = \frac{1}{r} w_r - \frac{1}{r^2} w, \qquad
r^2 u_r = r w_r - w \\[1ex]
\big(r^2 u_r \big)_r = r w_{rr}, \qquad 
\frac{1}{r^2} \big(r^2 u_r \big)_r = \frac{1}{r} w_{rr} 
\end{dcases}
\label{w_deriv}
\end{equation}
Substituting the expressions from \eqref{w_deriv} into \eqref{wav_eq3} yields 
the IBVP for $w(r, t)$. 
\begin{equation}
\begin{dcases}
w_{tt}(r, t) = c^2 w_{rr}(r, t),  
\quad r \in (0, \infty), \; t > 0\\[1ex]
w(0, t) = 0 \\
w(r, 0) = 0, \quad w_t(r, 0) = r \sigma(r)
\end{dcases}
\label{wav_eq4B}
\end{equation}
To satisfy the boundary condition $w(0, t) = 0$ in \eqref{wav_eq4B}, 
we extend both the unknown function $w(r, t)$ and the initial condition $g(r) \equiv r \sigma(r)$ 
to odd functions of spatial variable $r$. 
Let $g_\text{odd} (r) $ be the odd extension of $g(r)$. With this extension, the extended function 
$w(r, t)$ satisfies a pure initial value problem (IVP) on the 1D wave equation. 
\begin{equation}
\begin{dcases}
w_{tt}(r, t) = c^2 w_{rr}(r, t),  
\quad r \in (-\infty, +\infty), \; t > 0\\[1ex]
w(r, 0) = 0, \quad w_t(r, 0) = g_\text{odd}(r)
\end{dcases}
\label{wav_eq5}
\end{equation}
Let $G(s) \equiv \int_{+\infty}^{s} g_\text{odd}(s') d s'
= \int_{+\infty}^{|s|} s' \sigma(s') ds' $. $G(s) $ satisfies 
\begin{equation}
G'(s) = g_\text{odd}(s), \quad G(-\infty) = G(+\infty) = 0, \quad G(-s) = G(s)
\end{equation}
It is straightforward to verify that $w(r, t)$ given below satisfies IVP \eqref{wav_eq5}
and IBVP \eqref{wav_eq4B}:
\begin{equation}
\boxed{\quad 
w(r, t) = \frac{1}{2c}\Big(G(r+ct)-G(r-ct) \Big), \quad 
G(s) \equiv \int_{+\infty}^{|s|} s' \sigma(s') d s'
\quad }
\label{sol_eq5}
\end{equation}
This is known as d'Alembert's solution. 
We now go from solution $w(r, t)$ of IBVP \eqref{wav_eq4B} to solution $u^\text{(IV)}(r, t)$ 
of IBVP \eqref{wav_eq3} using the relation between the two given in \eqref{w_deriv}. 
\[ u^\text{(IV)}(r, t) = \frac{1}{r}w(r, t) \]
Substituting this expression into \eqref{Duhamel}, we obtain the solution of 
\eqref{wav_eq2}: 
\begin{equation}
\boxed{\quad 
\begin{dcases} 
u(r, t) = \frac{1}{2c\, r} \int_0^t f(t-\tau) 
\Big(G(r+c\tau)-G(r-c\tau) \Big) d\tau, \quad r > 0, \; t > 0 \\[1ex]
\qquad \text{Forcing: } F(r, t) = \sigma(r) f(t), \qquad 
G(s) \equiv \int_{+\infty}^{|s|} s' \sigma(s') d s' 
\end{dcases} \quad}
\label{eq2_sol}
\end{equation}
We remark that our derivation of the wave equation under spherical symmetry, along with the application of Duhamel's principle and the use of odd extensions, follows standard techniques as described in \cite{Evans_2010, John_1982, Strauss_2007}.

\subsection{Numerical Procedure for Evaluating \eqref{eq2_sol}} 
The solution of \eqref{wav_eq2} with force $\sigma(r) f(t)$ is given in \eqref{eq2_sol}.
It contains two layers of integration. For general functions $\sigma(r)$ and $f(t)$, evaluating these integrals requires numerical computation. 
Below we describe the numerical procedure for evaluating \eqref{eq2_sol}.

We consider the situation where the spatial profile of the force $\sigma(r)$
has compact support. Specifically, we assume that 
there exists $d > 0$ such that \(\sigma(r) = 0, \; r \in [d, +\infty) \). In other words, the active force is confined within a sphere of radius $d$ centered at the origin. Consequently, the function $G(s)$ satisfies \(G(s) = 0 \) for \( s \in [d, +\infty) \), 
and in particular, $G(d) = 0$. 
The method for evaluating \eqref{eq2_sol} consists of the algorithms listed below.
\vspace*{2ex}
\begin{algorithm}
\caption{Cubic spline representation of $G(s)$, for $ s \in [0, d]$} \label{alg_1}
\begin{algorithmic}[1]
\Function{get\_pp}{$d, b, N$} \Comment{$\sigma(\;)$ is a function.} 
\State $g \gets \text{lambda $s$:}\; s\,\sigma(s)$ \Comment{$g(\;)$ is a function.}
\State $h \gets d/N$ \Comment{$h $ is the spatial step in integration.}
\State $\{s_j\} \gets \{j h, \; 0 \le j \le N \}$
\State $\{g_j\} \gets g(\{s_j\})$
\State $\{g_{j-\frac{1}{2}}\} \gets g(\{s_j \}-\frac{1}{2} h)$
\Comment{$1 \le j \le N$}
\State $\{G_j\} \gets \big[0, \; \text{cumsum}\big(\{g_{j-1}+4g_{j-\frac{1}{2}}+g_j\}
\big) h/6 \big] $ \Comment{$1 \le j \le N$}
\State \Comment{$G_j \approx \int_0^{s_j} g(s') ds'$, obtained with Simpson's rule.}
\State $\{G_j\} \gets \{G_j\}-G_N$ \Comment{update $\{G_j\}$ to enforce $G(+\infty)=0$}
\State $pp \gets \text{spline}\big(\{s_j\}, \{G_j\}\big)$
\State \Return $pp$  \Comment{$pp$ is a cubic spline representation of $G(s)$.}
\EndFunction
\end{algorithmic}
\end{algorithm}

\begin{algorithm}
\caption{Construct the vectorized function $G(s)$ for all $s \in (-\infty, +\infty)$}\label{alg_2}
\begin{algorithmic}[1]
\Function{G}{$s, d, pp$} \Comment{$s = \{s_j\}$ is a vector.} 
\State $s_d \gets \min\big(d, | s | \big)$ \Comment{$s_d$ is a vector.}
\State $y \gets \text{ppval}(pp, s_d)$ \Comment{$y$ is a vector.}
\State \Return $y$  \Comment{$y \approx G(s)$ for any vector input $s$.}
\EndFunction
\end{algorithmic}
\end{algorithm}

\begin{algorithm}
\caption{Compute $u(r, T)$ for given scalar values of $r$ and $T$}\label{alg_3}
\begin{algorithmic}[1]
\Function{u}{$r, T, c, d, b, f, N, h_u$} 
\Comment{$c$ is sound speed; $\sigma(\; )$ and $f(\; )$ are functions.} 
\State \Comment{$h_u$ is the time step specified by user.} 
\State $t_b \gets \min\big(T, \max\big(0, (r-d)/c \big) \big)$ \Comment{$t_b$ is the lower limit of
$\tau$ in \eqref{eq2_sol}}
\State $t_e \gets \min\big(T, (r+d)/c \big)$ \Comment{$t_e$ is the upper limit of
$\tau$ in \eqref{eq2_sol}}
\State $N_t \gets \max\big(1, \text{ceil}((t_e-t_b)/h_u)\, \big) $ 
\Comment{$N_t \ge 1$ even if $t_e=t_b$} 
\State $h \gets (t_e-t_b)/N_t$ \Comment{$h$ is the time step used in integration}
\State $\{t_k\} \gets \{t_b+k h, \; 0 \le k \le N_t \}$
\State $\{t_{k-\frac{1}{2}}\} \gets \{t_k-\frac{1}{2}h, \; 1 \le k \le N_t \}$
\State 
\State $pp \gets \Call{get\_pp}{d, b, N}$
\State $\{w_k\} \gets f(T-\{t_k\}) \Big(\Call{G}{r+c\{t_k\}, d, pp}
-\Call{G}{r-c\{t_k\}, d, pp}\Big)$
\State $\{w_{k-\frac{1}{2}}\} \gets f(T-\{t_{k-\frac{1}{2}}\}) 
\Big( \Call{G}{r+c \{t_{k-\frac{1}{2}}\}, d, pp}
-\Call{G}{r-c \{t_{k-\frac{1}{2}}\}, d, pp} \Big) $
\State $y \gets \text{sum}\big(\{w_{k-1}+4w_{k-\frac{1}{2}}+w_k\}\big) 
(h/6)/(2c\, r) $ 
\Comment{$1 \le k \le N_t$}
\State \Return $y$  \Comment{$y \approx u(r, T)$ given in \eqref{eq2_sol}}
\EndFunction
\end{algorithmic}
\end{algorithm}
\clearpage

\subsection{Asymptotic Quasi-Steady Periodic Oscillations} \label{QS_oscillation}
In the present study, the force is given by $F(r, t) = \sigma(r) f(t)$. 
When the spatial forcing profile $\sigma(r)$ has compact support, the corresponding function $G(s)$, 
defined in \eqref{eq2_sol}, also has compact support. Specifically, 
$\sigma(r) = 0 \text{ for } r \ge d$ implies 
\begin{equation}
 G(s) = \int_{+\infty}^{|s|} s' \sigma(s') d s' = 0\; \text{ for } |s| \ge d
 \end{equation}
When the spatial forcing profile $\sigma(r)$ virtually has compact support, the corresponding 
$G(s)$ has the same behavior. 
For example, consider the Gaussian spatial profile $\sigma(r) = \frac{1}{d^3}\rho(\frac{r}{d})$ 
where $\rho(s)$ is the standard normal density. The support of this $\sigma(r)$
is effectively limited to a small multiple of $d$. The corresponding $G(s)$ is given by
\begin{equation}
\begin{dcases}
\sigma(s) = \frac{1}{d^3 \sqrt{2\pi}} e^{\frac{-s^2}{2d^2}}, \quad 
\int s \sigma(s)ds= -d^2 \sigma(s) + C \\[1ex]
G(s) = \int_{+\infty}^{|s|} s' \sigma(s') d s' = -d^2 \sigma(|s|) 
= \frac{-1}{d \sqrt{2\pi}} e^{\frac{-s^2}{2d^2}}
\end{dcases} \label{G_cal}
\end{equation}
In this example, $\sigma(r)$ and $G(s)$ have the same support. 
Generally, when the spatial profile $\sigma(r)$ is effectively confined to a finite spatial region 
$(r < R)$, the effective support of the corresponding $G(s)$ is $(-R, R)$, 
which we call the forcing core. 

We examine the solution given in \eqref{eq2_sol} at a spatial location $r$ outside the forcing 
core, at a sufficiently large time $t$ such that the effect of starting forcing at $t=0$ 
has propagated beyond the location $r$. In this situation, 
the solution exhibits the same behavior as if the force has started at $t = -\infty$. 
In the force expression $F(r, t) =\sigma(r)f(t)$, when the time profile $f(t)$ is periodic, 
the solution exhibits a quasi-steady oscillation at large times. 
The quasi-steady state is in the sense that the solution is approximately periodic. 
If the spatial profile has true compact support (i.e., being exactly zero outside the support), 
then the quasi-steady solution is exactly periodic. 

In \eqref{eq2_sol}, we decompose the general solution into 
the outward propagating part $u_\text{out}$ and the inward propagating part $u_\text{in}$. 
To facilitate analysis, we apply a change of variables in each component to center the integrand around $G(s)$. 
\begin{align}
& u(r, t) = \underbrace{(\frac{-1}{2c\, r})\int_0^t f(t-\tau) G(r-c\tau) d\tau
}_{u_\text{out}} 
+\underbrace{(\frac{1}{2c\, r}) \int_0^t f(t-\tau) G(r+c\tau) d\tau
}_{u_\text{in}} \nonumber \\
& \quad= \underbrace{(\frac{-1}{2c^2 r})\int_{-(ct -r)}^r f(t-\frac{r}{c}+\frac{s'}{c}) G(s') ds'
}_{u_\text{out}} 
+\underbrace{(\frac{1}{2c^2r}) \int_r^{r+ct} f(t+\frac{r}{c}-\frac{s'}{c}) G(s') ds'
}_{u_\text{in}} \label{u_out_in}
\end{align}
Mathematically, we analyze the solution $u(r, t)$ in the quasi-steady state region of $(r, t)$, 
denoted by $D_{(r, t)}^\text{(QSS)}$ and defined as: 
\begin{equation} D_{(r, t)}^\text{(QSS)} \equiv 
\Big\{(r, t) \Big| r > R, \; (-R + c t ) >r \Big\}
\label{r_t_region}
\end{equation} 
In region $D_{(r, t)}^\text{(QSS)}$, $r > R$, and as a result, the integration interval in the inward propagating component $u_\text{in}$ lies entirely outside the effective support of $G(s)$. 
Because $G(s)$ is very small outside the effective support (for example, consider the Gaussian function outside 5 times the standard deviation), the contribution from $u_\text{in}$ is negligible: 
\begin{equation}
\begin{dcases}
 [r, \, r+ct] \, \cap \, [-R, R] = \text{null}, \\[1ex]
 u_\text{in}=(\frac{1}{2c^2r}) \int_r^{r+ct} f(t+\frac{r}{c}-\frac{s'}{c}) G(s') ds'
\approx 0
\end{dcases} \label{u_in_asy}
\end{equation}
In region $D_{(r, t)}^\text{(QSS)}$, $r > R$ and $(-R + c t ) >r$, leading to 
$[-(ct-r), r] \supset [-R, R] $. That is, the integration interval 
in the outward propagating component $u_\text{out}$ covers the effective support 
of $G(s)$. Thus, we can approximate it using the integral over $(-\infty, +\infty)$. 
\begin{equation}
\begin{dcases}
[-(ct-r), r] \supset [-R, R] , \\[1ex]
u_\text{out}
\approx (\frac{-1}{2c^2 r})\int_{-\infty}^{+\infty} f(t-\frac{r}{c}+\frac{s'}{c}) G(s') ds'
\end{dcases} \label{u_out_asy}
\end{equation}
We consider a sinusoidal oscillating force, $f(t) = \sin(\omega t)$. 
Substituting this into \eqref{u_out_asy} and utilizing the fact that the function $G(s)$ is even, we write $u_\text{out}$ as 
\begin{align*}
u_\text{out} & \approx 
(\frac{-1}{2c^2 r})\int_{-\infty}^{+\infty} \sin\big(\omega (t-\frac{r}{c})+
\frac{\omega}{c}s' \big) G(s') ds' \\
& = (\frac{-1}{2c^2}) \frac{1}{r} \sin\big(\omega (t-\frac{r}{c})\big) 
\int_{-\infty}^{+\infty} \cos\big(\frac{\omega}{c}s' \big) G(s') ds'
\end{align*}
In summary, for $(x, t) \in D_{(r,t)}^\text{(QSS)}$ defined in \eqref{r_t_region}, 
the integral solution of the 3D forced wave equation given in \eqref{eq2_sol}
is asymptotically a quasi-steady periodic oscillation: 
\begin{equation}
\boxed{\quad 
\begin{dcases} 
u(r, t) \approx (\frac{-1}{2c^2}) \frac{1}{r} \sin\big(\omega (t-\frac{r}{c})\big) 
\int_{-\infty}^{+\infty} \cos\big(\frac{\omega}{c}s' \big) G(s') ds' \\[1ex]
\qquad \text{ for } (r, t) \in D_{(r, t)}^\text{(QSS)} \equiv 
\Big\{(r, t) \Big| r > R, \; (-R + c t ) >r \Big\}, \\[1ex]
\qquad \text{ and } F(r, t) = \sigma(r) \sin(\omega t).
\end{dcases} \quad}
\label{eq2_sol_asy}
\end{equation}

\subsection{Exact Solution of Case 1}
\label{Exact_sol}
\begin{figure}[!h]
\vskip 0.5cm
\begin{center}
\psfig{figure=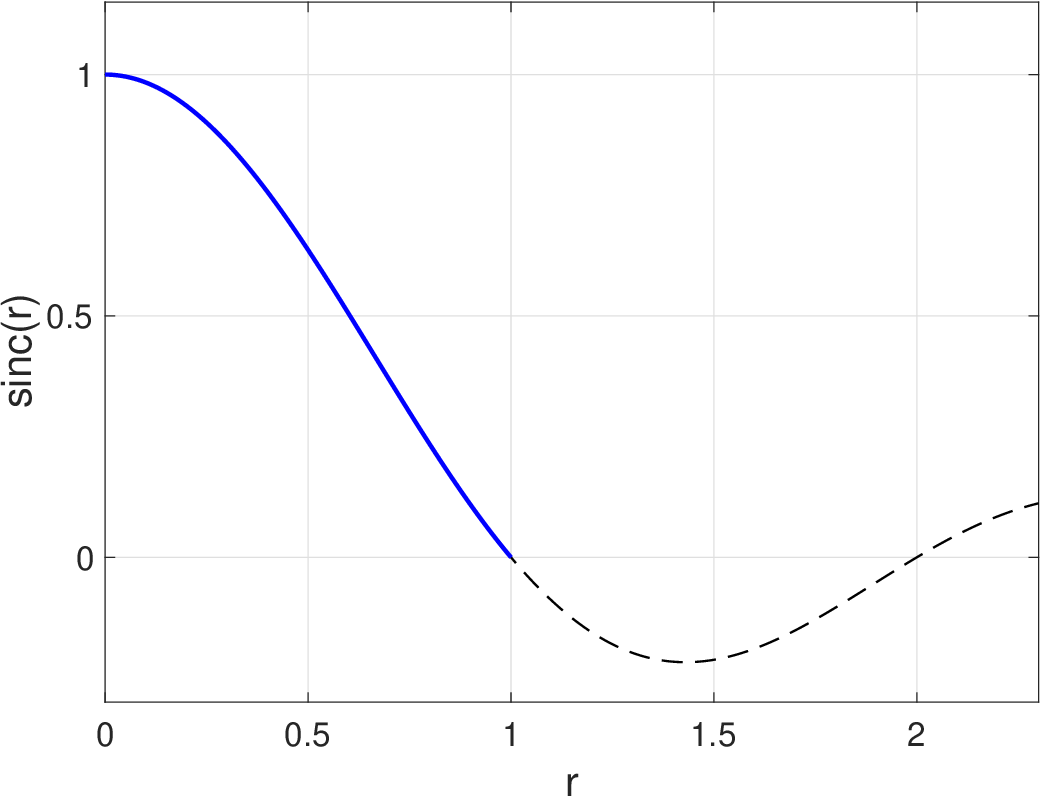, width=3.5in}
\end{center}
\vskip -0.75cm
\caption{Graphs of the normalized sinc function $\text{sinc}(r) $ and its truncated version.}
\label{fig_01}
\end{figure}

We consider the special case where the forcing spatial profile $\sigma(r)$ is 
a truncated and scaled sinc function. Mathematically, the force has the expression:
\begin{equation}
\begin{dcases} 
F(r, t) = \sigma(r) \sin(\omega t) \\
\sigma(r) = \frac{1}{d^3} \text{sinc}^{(tc)}(\frac{r}{d}), \qquad 
\text{sinc}^{(tc)}(s) \equiv \begin{dcases}
\text{sinc}(s), & 0 \le s \le 1 \\
0, & s \ge 1
\end{dcases} 
\end{dcases} 
\label{f_case1} 
\end{equation}
where $\text{sinc}(s) \equiv \frac{\sin(\pi s)}{\pi s}$ is the normalized sinc 
function. Graphs of $\text{sinc}(r) $ and its truncated version are shown in Figure \ref{fig_01}. 
In this special case, we evaluate the integrals in \eqref{eq2_sol} 
analytically, without any approximation. 
For the given $\sigma(r)$, we have 
\begin{align}
& G(s) = \int_d^{|s|} s' \sigma(s') ds' = \frac{1}{\pi^2 d} 
\int_d^{|s|} \sin(\frac{\pi}{d} s') d(\frac{\pi}{d}s') , \quad s \in [-d, d] \nonumber \\[1ex]
& \qquad = \frac{-1}{\pi^2 d} \Big(1+\cos(\frac{\pi}{d}s) \Big) 
\label{G_exact} 
\end{align}
The exact support of $G(s)$ is $[-d, d]$. 
In \eqref{u_out_in}, we identified the two integrals as the outward-propagating 
and inward-propagating. Now, using the exact supports of 
$G(r-c\tau)$ and $G(r+c\tau)$, we determine the precise integration intervals for these two integrals. 
\begin{equation}
\begin{dcases}
u(r, t) = \underbrace{(\frac{-1}{2c\, r})\int_{t_b}^{t_e} f(t-\tau) G(r-c\tau) d\tau
}_{u_\text{out}} 
+\underbrace{(\frac{1}{2c\, r}) \int_{t_1}^{t_2} f(t-\tau) G(r+c\tau) d\tau
}_{u_\text{in}} \\
 t_b = \min(t, \max(0, (r-d)/c)), \qquad t_e = \min(t, (r+d)/c) \\
 t_1 = \min(t, \max(0, -(r+d)/c))=0, \qquad t_2 = \min(t, \max(0, -(r-d)/c)) 
 \end{dcases} \label{u_int_limits}
\end{equation}
Using $G(s)$ from \eqref{G_exact} and $f(t) = \sin(\omega t)$, we calculate
$u_\text{in}$ and $u_\text{out}$ analytically 
\begin{align}
& u_\text{out}= \frac{1}{2\pi^2 c d}\, \frac{1}{r} \int_{t_b}^{t_e}
\sin(\omega(t-\tau)) \Big(1+\cos(\frac{\pi}{d}(r-c\tau)) \Big)d\tau \nonumber \\
& \quad =\frac{1}{2\pi^2 c d}\, \frac{1}{r} \Big[ 
\int_{t_b}^{t_e}\sin(\omega t-\omega \tau) d\tau
+\frac{1}{2}\int_{t_b}^{t_e} \sin((\omega t + \frac{\pi r}{d})
-(\omega+\frac{\pi c}{d})\tau) d\tau \nonumber \\
& \hspace{1in} +\frac{1}{2}\int_{t_b}^{t_e} \sin((\omega t - \frac{\pi r}{d})
-(\omega-\frac{\pi c}{d})\tau) d\tau \Big] \label{int_res1} \\
& u_\text{in}  = 
\frac{-1}{2\pi^2 c d}\, \frac{1}{r}\int_{t_1}^{t_2} \sin(\omega(t-\tau)) 
\Big(1+\cos(\frac{\pi}{d}(r+c\tau)) \Big)d\tau \nonumber \\
& \quad =\frac{-1}{2\pi^2 c d}\, \frac{1}{r}\Big[ 
 \int_{t_1}^{t_2}\sin(\omega t-\omega \tau) d\tau
+\frac{1}{2}\int_{t_1}^{t_2} \sin((\omega t + \frac{\pi r}{d})
-(\omega-\frac{\pi c}{d})\tau) d\tau \nonumber \\
& \hspace{1in} +\frac{1}{2}\int_{t_1}^{t_2} \sin((\omega t - \frac{\pi r}{d})
-(\omega+\frac{\pi c}{d})\tau) d\tau \Big] \label{int_res2} 
\end{align}
All integrals in \eqref{int_res1} are of the form
$\int_{t_b}^{t_e} \sin(a_0-v \tau) d\tau $, for which we derive the formula below:
\begin{align}
& \int_{t_b}^{t_e} \sin(a_0-v \tau) d\tau 
= \frac{1}{v}\big(\cos(a_0-v t_e)-\cos(a_0-v t_b)\big) \nonumber \\
& \qquad = 2 \sin(a_0-v\frac{(t_b+t_e)}{2})
\frac{\sin(v\frac{(t_e-t_b)}{2})}{v} \nonumber \\[1ex]
& \qquad =(t_e-t_b) \text{sinc}(\frac{v}{\pi}\frac{(t_e-t_b)}{2}) 
\sin(a_0-v\frac{(t_b+t_e)}{2}) \label{int_formu}  
\end{align}
We apply formula \eqref{int_formu} to the integrals in \eqref{int_res1} to evaluate 
$u_\text{out}$ and calculate $u_\text{in}$ in (\ref{int_res2}) in a similar manner. 
Solution $u(r, t)$ is the sum of $u_\text{out}$ and $u_\text{in}$.
\begin{equation}
\boxed{\quad \begin{aligned}
& u_\text{out} =\frac{(t_e-t_b)}{2\pi^2 c d}\, \frac{1}{r} \Big[
\text{sinc}(\frac{\omega}{\pi}\frac{(t_e-t_b)}{2})
\sin\big(\omega (t-\frac{(t_b+t_e)}{2}) \big)
\\[1ex]
& \qquad + \frac{1}{2}\text{sinc}(\frac{(\omega +\frac{\pi c}{d})}{\pi}\frac{(t_e-t_b)}{2})
\sin\big (\omega (t-\frac{(t_b+t_e)}{2}) + \frac{\pi}{d} (r-c\frac{(t_b+t_e)}{2}) \big) 
\\[1ex]
& \qquad +\frac{1}{2}\text{sinc}(\frac{(\omega -\frac{\pi c}{d})}{\pi}\frac{(t_e-t_b)}{2}) 
\sin\big (\omega (t-\frac{(t_b+t_e)}{2}) - \frac{\pi}{d} (r-c\frac{(t_b+t_e)}{2}) \big) 
\Big] \\[1ex]
& \quad t_b = \min(t, \max(0, (r-d)/c)), \qquad t_e = \min(t, (r+d)/c) 
\end{aligned} \quad}\label{u_out_exp}
\end{equation}
\vskip 1ex
\begin{equation}
\boxed{\quad \begin{aligned}
& u_\text{in} = \frac{-t_2}{2\pi^2 c d}\, \frac{1}{r}\Big[
\text{sinc}(\frac{\omega}{\pi}\frac{t_2}{2})
\sin\big(\omega (t-\frac{t_2}{2}) \big)
\\[1ex]
& \qquad + \frac{1}{2}\text{sinc}(\frac{(\omega -\frac{\pi c}{d})}{\pi}\frac{t_2}{2}) 
\sin\big(\omega (t-\frac{t_2}{2}) + \frac{\pi}{d} (r+c\frac{t_2}{2}) \big) 
\\[1ex]
& \qquad +\frac{1}{2}\text{sinc}(\frac{(\omega +\frac{\pi c}{d})}{\pi}\frac{t_2}{2}) 
\sin\big(\omega (t-\frac{t_2}{2}) - \frac{\pi}{d} (r+c\frac{t_2}{2}) \big) 
\Big] \\[1ex]
& \quad t_2 = \min(t, \max(0, -(r-d)/c)) 
\end{aligned} \quad} \label{u_in_exp}
\end{equation}
\begin{equation}
\boxed{\qquad \begin{dcases}  u(r, t) = u_\text{out} + u_\text{in} 
\;\; \text{ for } r \ge 0, \; t \ge 0, \; \text{ and } \\[1ex]
\qquad F(r, t) = \frac{1}{d^3} \text{sinc}(\frac{r}{d}) 
\sin(\omega t), \;\; r \le d. 
\end{dcases}
\qquad} \label{u_special}
\end{equation}
\vskip 1ex 
The analytical solution of the 3D forced wave equation in case 1 is given in 
\eqref{u_special} with \eqref{u_out_exp} and \eqref{u_in_exp}; 
it is exactly valid for all $r > 0$ and $t > 0$; 
there is no approximation involved; the effect of starting forcing at $t=0$ is 
contained in $(t_b, t_e, t_2)$. 

\subsection{Quasi-Steady Periodic Solutions in Several Cases}
\label{asy_sol}
In this subsection, we calculate the asymptotic quasi-steady periodic solution 
given in \eqref{eq2_sol_asy} with time profile $f(t) = \sin(\omega t)$ 
for several cases of spatial profile $\sigma(r)$. 

\subsubsection{Case 1: $\sigma(r)$ is a truncated $\text{sinc}(\;)$ function}
\label{case_1}
The force of case 1 is given in \eqref{f_case1}.
In \eqref{u_special} with \eqref{u_out_exp} and \eqref{u_in_exp}, 
we derived the exact solution of case 1, which is exactly 
valid for all $r > 0$ and $t > 0$, 
at the price that $(t_b, t_e, t_2)$ has complicated expressions. 
It is not straightforward to decipher the behavior of the exact solution at a general point $(r, t)$. 
In contrast, the general asymptotic solution given in \eqref{eq2_sol_asy} exhibits a clear 
periodic oscillation. 
We now derive the asymptotic solution of case 1 from the exact solutions 
\eqref{u_out_exp} and \eqref{u_in_exp}. 
We show that in case 1, the asymptotic solution is exactly valid in region 
$D_{(r, t)}^\text{(QSS)}$. 

In case 1, the support of $G(s)$ is $(-d, d)$, and the $(r, t)$-region 
$D_{(r, t)}^\text{(QSS)}$, introduced in \eqref{r_t_region} when deriving the general 
asymptotic solution, becomes  
\begin{equation}
 D_{(r, t)}^\text{(QSS)} \equiv 
\Big\{(r, t) \Big| r > d, \; (-d + c t ) >r \Big\} \label{r_t_region_cs1}
\end{equation}
For $(r, t) \in D_{(r, t)}^\text{(QSS)}$, the expressions
in \eqref{u_out_exp} and \eqref{u_in_exp} yield: 
\begin{equation}
\begin{dcases}
t_2 = 0, \qquad u_\text{in} = 0 \\[1ex]
 t_b = \frac{r-d}{c}, \qquad 
t_e = \frac{r+d}{c} \\[1ex]
t-\frac{t_b+t_e}{2} = t-\frac{r}{c}, \qquad 
r-c\frac{t_b+t_e}{2} = 0, \qquad 
\frac{t_e-t_b}{2} = \frac{d}{c} \\[1ex]
u_\text{out} =  \frac{1}{\pi^2 c^2}\Big(
\text{sinc}(\frac{\omega d}{\pi c}) 
+\frac{1}{2}\text{sinc}(\frac{\omega d}{\pi c}+1)
+\frac{1}{2}\text{sinc}(\frac{\omega d}{\pi c}-1) \Big) 
\frac{1}{r}\sin(\omega (t- \frac{r}{c})) 
\end{dcases}
\label{t_limits}
\end{equation}
Note that for $(r, t) \in D_{(r, t)}^\text{(QSS)}$, the expressions in 
\eqref{t_limits} are exactly valid with no approximation. 
We write the solution in terms of $\omega$ and $\lambda$ respectively. 
\begin{equation}
\boxed{\quad \begin{dcases} 
\text{Solution in terms of angular frequency $\omega$:} \\[1ex]
u(r, t) = \underbrace{\frac{1}{\pi^2 c^2}\Big(
\text{sinc}(\frac{\omega d}{\pi c}) 
+\frac{1}{2}\text{sinc}(\frac{\omega d}{\pi c}+1)
+\frac{1}{2}\text{sinc}(\frac{\omega d}{\pi c}-1) \Big)}_{\text{independent of $(r, t)$}}
\frac{1}{r}\sin(\omega (t- \frac{r}{c})) \\[1ex]
\text{Solution in terms of wavelength 
$\displaystyle \lambda = \frac{2\pi c}{\omega}$:} \\[1ex] 
u(r, t) = \underbrace{\frac{1}{\pi^2 c^2}\Big(
\text{sinc}(\frac{2 d}{\lambda}) 
+\frac{1}{2}\text{sinc}(\frac{2 d}{\lambda}+1)
+\frac{1}{2}\text{sinc}(\frac{2 d}{\lambda}-1) \Big)}_{\text{independent of $(r, t)$}}
\frac{1}{r}\sin(2\pi \frac{(c t- r)}{\lambda}) \\[1ex]
\qquad \text{ for } (r, t) \in D_{(r, t)}^\text{(QSS)} \equiv 
\Big\{(r, t) \Big| r > d, \; (-d + c t ) >r \Big\}, \\[1ex]
\qquad \text{ and } F(r, t) = \frac{1}{d^3} 
\text{sinc}(\frac{r}{d}) \sin(\omega t), \;\; r \le d.
\end{dcases} \quad }
 \label{u_sinc_QSS}
\end{equation}
In region $D_{(r, t)}^\text{(QSS)}$, the asymptotic solution of case 1 given in 
\eqref{u_sinc_QSS} is the same as the exact solution. There is no approximation involved. 

\subsubsection{Case 2: $\sigma(r)$ is a Gaussian function} 
\label{case_2}
In case 2, the force has the expression: 
\begin{equation}
\begin{dcases} 
F(r, t) = \sigma(r) \sin(\omega t) \\
\sigma(r) = \frac{1}{d^3}\rho(\frac{r}{d}), \qquad 
\rho(s) \equiv \frac{1}{\sqrt{2\pi}} e^{\frac{-s^2}{2}} 
\end{dcases} 
\label{f_case2} 
\end{equation}
Here, $\rho(s)$ is the standard normal density function, which satisfies: 
\[ \begin{dcases}
\rho(-s) = \rho(s), \qquad \int_{-\infty}^{+\infty} \rho(s) ds = 1
, \qquad \int_{-\infty}^{+\infty} \rho(s) s^2 ds = 1  \\
\rho'(s) = -s \rho(s), \qquad \int s \rho(s)ds = -\rho(s) + C 
\end{dcases} \]
We evaluate the asymptotic solution given in \eqref{eq2_sol_asy} 
for $\sigma(r)$ given in \eqref{f_case2}. Function $G(s)$ is 
\begin{align}
& G(s) = \int_{+\infty}^{s} s' \sigma(s') ds' 
= \frac{1}{d} \int_{+\infty}^{s} (\frac{s'}{d}) \rho(\frac{s'}{d}) d(\frac{s'}{d}) 
= -\frac{1}{d}\rho(\frac{s}{d}), \quad s \in (-\infty, +\infty)
\label{G_exact_2} 
\end{align}
The effective support of $G(s)$ is the interval $(-R, R)$ where $R/d \gg 1$. 
For example, with a relative error tolerance of $10^{-14}$, 
$G(s)$ can be considered effectively confined to $(-R, R)$ for $R = 8d$. 
The $(r, t)$-region for the asymptotic quasi-steady solution is given in \eqref{r_t_region}. 
To evaluate the integral in \eqref{eq2_sol_asy}, we derive a formula: 
\begin{align}
& \int_{-\infty}^{+\infty} \cos(v s) \rho(s) ds
= \text{Re}\Big( \int_{-\infty}^{+\infty} 
\frac{1}{\sqrt{2\pi}}e^{\frac{-1}{2}s^2+i v s} ds \Big) \nonumber \\
& \quad 
= \text{Re}\Big( e^{-\frac{1}{2}v^2} \underbrace{\int_{-\infty}^{+\infty} 
\frac{1}{\sqrt{2\pi}}e^{\frac{-1}{2}(s-i v)^2} ds}_{=1} \Big) 
= e^{-\frac{1}{2}v^2} \label{int_formu2}
\end{align}
We apply formula \eqref{int_formu2} to the integral in \eqref{eq2_sol_asy} to derive 
asymptotic solution $u(r, t)$. 
\begin{align*} 
& \int_{-\infty}^{+\infty} \cos\big(\frac{\omega}{c}s' \big) G(s') ds' 
= -\int_{-\infty}^{+\infty} \cos\big(\frac{\omega d}{c}s \big) \rho(s) ds
=- e^{\frac{-\omega^2 d^2}{2c^2}} \\[1ex]
& u(r, t) \approx \frac{1}{2c^2}e^{\frac{-\omega^2 d^2}{2c^2}} 
\frac{1}{r} \sin\big(\omega (t-\frac{r}{c})\big) 
\end{align*} 
We write the solution respectively in terms of $\omega$ and $\lambda$. 
\begin{equation}
\boxed{\quad 
\begin{dcases}
\text{Solution in terms of angular frequency $\omega$:} \\[1ex] 
\quad u(r, t) \approx \frac{1}{2c^2}e^{\frac{-\omega^2 d^2}{2c^2}}
\frac{1}{r}\sin(\omega (t- \frac{r}{c})) \\[1ex]
\text{Solution in terms of wavelength 
$\displaystyle \lambda = \frac{2\pi c}{\omega}$:} \\[1ex] 
\quad u(r, t) \approx \frac{1}{2c^2}e^{-2\pi^2\frac{d^2}{\lambda^2}}
\frac{1}{r}\sin(2\pi \frac{(c t- r)}{\lambda}) \\[1ex]
\qquad \text{ for } (r, t) \in D_{(r, t)}^\text{(QSS)} \equiv 
\Big\{(r, t) \Big| r > R, \; (-R + c t ) >r \Big\}, \;\; R/d \gg 1 , \\[1ex]
\qquad \text{ and } F(r, t) = \frac{1}{d^3} \rho(\frac{r}{d}) \sin(\omega t).
\end{dcases} \quad} 
\label{u_normal_QSS}
\end{equation}
Solution of case 2 given in \eqref{u_normal_QSS} is 
approximately valid in region $D_{(r, t)}^\text{(QSS)}$ with 
$R/d \gg 1$. More specifically, the error of approximation decreases exponentially as
the ratio $r/d$ increases, similar to the decay of the tail in a Gaussian function.

\subsubsection{Case 3: $\sigma(\vec{x})$ is a Dirac delta function (point-force)}
\label{case_3}
In case 3, the force has the expression 
\begin{equation}
F(r, t) = \sigma(r) \sin(\omega t), \quad  
\sigma(r) = \delta(\vec{x})  
\label{f_case3} 
\end{equation}
In \eqref{f_case3}, $\delta(\vec{x})$ needs to be interpreted as the delta function of 
$\vec{x} \in \mathbb{R}^3$. Although $\delta(\vec{x})$ is spherically symmetric, 
we shall not write it as $\delta(r)$, which is a completely different function. 
Case 3 is related to the limit of case 2. As $d \rightarrow 0_+$, 
spatial profile $\sigma_\text{case-2}(r) = \frac{1}{d^3}\rho(\frac{r}{d})$ in case 2 
converges to a multiple of the delta function $\delta(\vec{x})$. 
The multiplier is determined by the integral of $\sigma_\text{case-2}(|\vec{x}|)$.
\[ \int_{\mathbb{R}^3} \sigma_\text{case-2}(|\vec{x}|) d\vec{x} 
= \frac{1}{d^3} \int_0^{+\infty} \rho(\frac{r}{d}) 4 \pi r^2 dr
=  4\pi \int_0^{+\infty} \rho(s) s^2 ds = 2\pi \]
The spatial profiles of case 3 and case 2 are related by 
\begin{equation}
\underbrace{\; \sigma(\vec{x}) \;}_\text{case 3} = \frac{1}{2\pi} \lim_{d \rightarrow 0^+} 
\sigma_\text{case-2}(\vec{x})
\label{sgm_c2_c3}
\end{equation} 
It follows that the solutions of case 3 and case 2 are related by 
\begin{equation}
\underbrace{\; u(r, t) \;}_\text{case 3} = \frac{1}{2\pi} \lim_{d \rightarrow 0^+} 
u_\text{case-2}(r, t)
\label{u_c2_c3}
\end{equation} 
As $d \rightarrow 0_+$, the solution of case 2 given in \eqref{u_normal_QSS} converges to 
\begin{equation}
\lim_{d \rightarrow 0^+} u_\text{case-2}(r, t) =\frac{1}{2 c^2}
\frac{1}{r}\sin(\omega (t- \frac{r}{c})) 
\label{u_normal_lim}
\end{equation}
Combining \eqref{u_c2_c3} and \eqref{u_normal_lim}, 
we write the solution in terms of $\omega$ and $\lambda$ respectively. 
\begin{equation}
\boxed{\quad 
\begin{dcases} 
\text{Solution in terms of angular frequency $\omega$:} \\[1ex] 
\quad u(r, t) =\frac{1}{4 \pi c^2} \frac{1}{r}\sin(\omega (t- \frac{r}{c})) \\[1ex]
\text{Solution in terms of wavelength 
$\displaystyle \lambda = \frac{2\pi c}{\omega}$:} \\[1ex] 
\quad u(r, t) =\frac{1}{4 \pi c^2} \frac{1}{r}\sin(2\pi \frac{(c t- r)}{\lambda}) \\[1ex]
\qquad \text{ for } (r, t) \in D_{(r, t)}^\text{(QSS)} \equiv 
\Big\{(r, t) \Big| r > 0, \; c t > r \Big\}, \\[1ex]
\qquad \text{ and } F(\vec{x}, t) = \delta(\vec{x}) \sin(\omega t) 
\end{dcases} \quad} 
\label{u_delta_QSS}
\end{equation}
Solution of case 3 given in \eqref{u_delta_QSS} is exactly valid 
in region $D_{(r, t)}^\text{(QSS)}$, which is described by $c t > r $. 
Recall that in case 2, the error decays exponentially as the ratio $r/d$ increases. 
In case 3, $d = 0_+$ and $r/d = +\infty$ for any $r > 0$, leading to zero error in approximation. 

\subsubsection{Case 4: $\sigma(r) $ is concentrated in a thin spherical shell}
\label{case_4}
In case we just discussed, the spatial forcing is concentrated at one point. Now we consider 
the situation where the spatial forcing is concentrated in a spherical shell of infinitesimal thickness. 
In case 4, the force has the expression: 
\begin{equation}
F(r, t) = \sigma(r) \sin(\omega t), \quad \sigma(r) = \delta(r-r_0), \;\; r_0 > 0 
\label{f_case4} 
\end{equation}
The corresponding function $G(s)$ is 
\begin{align}
& G(s) = \int_{+\infty}^{|s|} s' \sigma(s') ds' 
= -\int_{|s|}^{+\infty} s' \delta(s'-r_0) ds' 
= -r_0, \quad s \in (-r_0, r_0)
\label{G_exact_4} 
\end{align}
The exact support of $G(s)$ is $(-r_0, r_0)$ and the $(r, t)$-region for 
the quasi-steady periodic solution, introduced in \eqref{r_t_region}, becomes 
\begin{equation}
 D_{(r, t)}^\text{(QSS)} \equiv 
\Big\{(r, t) \Big| r > r_0, \; (-r_0 + c t ) >r \Big\} \label{r_t_region_cs4}
\end{equation}
With $G(s) = \mathbf{1}_{[-r_0, r_0]}(s) $, 
the integral in general solution \eqref{eq2_sol_asy} becomes  
\begin{align*}
& \int_{-\infty}^{+\infty} \cos\big(\frac{\omega}{c}s' \big) G(s') ds' 
= -r_0\int_{-r_0}^{r_0} \cos\big(\frac{\omega}{c}s' \big) ds' 
= -2r_0\frac{c}{\omega} \sin(\frac{\omega }{c}r_0) 
\end{align*}
We write the solution in terms of $\omega$ and $\lambda$ respectively. 
\begin{equation}
\boxed{\quad 
\begin{dcases} 
\text{Solution in terms of angular frequency $\omega$:} \\[1ex] 
\quad u(r, t) = \frac{r_0}{\omega c} \sin(\frac{\omega }{c}r_0) 
\frac{1}{r} \sin\big(\omega (t-\frac{r}{c})\big) \\[1ex]
\text{Solution in terms of wavelength 
$\displaystyle \lambda = \frac{2\pi c}{\omega}$:} \\[1ex] 
\quad u(r, t) = \frac{r_0 \lambda}{2\pi c^2} \sin(2\pi \frac{r_0}{\lambda}) 
\frac{1}{r}\sin(2\pi \frac{(c t- r)}{\lambda}) \\[1ex]
\qquad \text{ for } (r, t) \in D_{(r, t)}^\text{(QSS)} \equiv 
\Big\{(r, t) \Big| r > r_0, \; (-r_0 + c t ) >r \Big\}, \\[1ex]
\qquad \text{ and } F(r, t) =\delta(r-r_0) \sin(\omega t), \quad r_0 > 0
\end{dcases} \quad} 
\label{u_del_r0_QSS}
\end{equation}
Similar to the situation of case 1 and case 3, the asymptotic solution of case 4 given in 
\eqref{u_del_r0_QSS} is the same as the exact solution in region $D_{(r, t)}^\text{(QSS)}$. 
This conclusion follows directly from the fact the compact support of function $G(s)$ is exact. 

Before we study masking, we summarize and clarify all solutions we obtained. 
\begin{itemize}
\item \eqref{eq2_sol} is the general exact solution of the 3D forced wave equation with 
$F(r, t) = \sigma(r) f(t)$ for arbitrary spatial profile $\sigma(r) $ and 
arbitrary time profile $f(t)$. Solution \eqref{eq2_sol} is expressed in two layers of integration 
and is valid for all $r > 0$ and $t > 0$. 
\eqref{eq2_sol} serves as the foundation for deriving other solutions.
\item \eqref{eq2_sol_asy} is the general asymptotic solution of the 3D forced wave equation
with $F(r, t) = \sigma(r) \sin(\omega t)$ for arbitrary spatial profile $\sigma(r) $. 
\eqref{eq2_sol_asy} is approximately valid in $D_{(r, t)}^\text{(QSS)}$, which describes 
the $(r, t)$-region in which position $r$ is outside the effective forcing core and time $t$ 
is large enough such that the effect of starting forcing at $t=0$ has propagated beyond position $r$. 
\eqref{eq2_sol_asy} is derived from the general exact solution \eqref{eq2_sol} and is used to derive asymptotic solutions for various spatial profiles $\sigma(r)$. 
\item \eqref{u_special} with \eqref{u_out_exp} and \eqref{u_in_exp} is the exact solution
for $ F(r, t) = \frac{1}{d^3} \text{sinc}(\frac{r}{d}) \sin(\omega t), \; r \le d$. 
It is a special case of \eqref{eq2_sol}. The specific expression of the force allows us to evaluate the integrals and write the result in terms of trigonometric functions. 
\item \eqref{u_sinc_QSS} is the asymptotic solution for 
$F(r, t) = \frac{1}{d^3} \text{sinc}(\frac{r}{d}) \sin(\omega t), \; r \le d$. 
It coincides with the exact solution in region $D_{(r, t)}^\text{(QSS)}$. 
\item \eqref{u_normal_QSS} is the asymptotic solution for 
$F(r, t) = \frac{1}{d^3} \rho_{N(0, 1)}(\frac{r}{d}) \sin(\omega t)$. 
The error of approximation decreases exponentially as
the ratio $r/d$ increases. 
\item \eqref{u_delta_QSS} is the asymptotic solution for 
$F(\vec{x}, t) = \delta(\vec{x}) \sin(\omega t)$. 
It coincides with the exact solution in region $D_{(r, t)}^\text{(QSS)}$, 
which is described by $c t > r $. 
\item \eqref{u_del_r0_QSS} is the asymptotic solution for 
$F(r, t) =\delta(r-r_0) \sin(\omega t), \; r_0 > 0$. 
It coincides with the exact solution in region $D_{(r, t)}^\text{(QSS)}$, 
which is described by $(-r_0 + c t) > r $. 
\end{itemize}
%

\section{Simple Examples of Masking}\label{masking_1}
This section illustrates simple examples of masking, highlighting both self-masking and masking 
by applying external point-force(s). 
\subsection{Self-masking in Case 1}
In \eqref{u_sinc_QSS}, we derived the asymptotic solution 
for the force $F(r, t) = \frac{1}{d^3}\text{sinc}(\frac{r}{d}) \sin(\omega t)$, $r \le d$. 
In region $D_{(r, t)}^\text{(QSS)}$, this asymptotic solution coincides with the exact solution. 
In \eqref{u_sinc_QSS}, when $d = \lambda \frac{k+1}{2}, \; k=1, 2, \ldots$, 
we have $\frac{2d}{\lambda} = k+1, \; k=1, 2, \ldots $ and 
\begin{equation}
\begin{dcases}
\text{sinc}(\frac{2d}{\lambda}) = \text{sinc}(\frac{2d}{\lambda}+1) 
= \text{sinc}(\frac{2d}{\lambda}-1) = \text{sinc}(k) = 0, \\[1ex]
u(r, t) \equiv 0 \quad 
\text{for all } (r, t) \in D_{(r, t)}^\text{(QSS)} \equiv 
\Big\{(r, t) \Big| r > d, \; (-d + c t ) >r \Big\} 
\end{dcases} 
\label{self_mask_c1}
\end{equation}
That is, for the force of case 1 with $d = \lambda \frac{k+1}{2}$ and $k=$ positive integer, 
the solution is identically zero outside the forcing core. 
This result exemplifies the phenomenon of \emph{self-masking} 
in which the wave field generated in the forcing core cancels itself in the quasi-steady region. 
Physically, this occurs due to destructive interference: at position $r$, the waves excited 
by different parts of the forcing core arrive with different phases; for a forcing core of
certain distribution at certain spatial scale relative to the wavelength, 
these waves exactly cancel each other,  
leading to that no wave outside the forcing core. 

\subsection{Self-masking in Case 4}
In \eqref{u_del_r0_QSS}, we derived the asymptotic solution 
for the force $F(r, t) =\delta(r-r_0) \sin(\omega t)$. 
In region $D_{(r, t)}^\text{(QSS)}$, this asymptotic solution coincides with the exact solution. 
In \eqref{u_del_r0_QSS}, when  
$r_0 = \lambda \frac{k}{2}, \; k=1, 2, \ldots$,  
we have $2\pi \frac{r_0}{\lambda} = k \pi, \; k=1, 2, \ldots$ and 
\begin{equation}
\begin{dcases}
\sin(2\pi \frac{r_0}{\lambda}) = \sin(k \pi) = 0, \\[1ex]
u(r, t) \equiv 0 \quad 
\text{for all } (r, t) \in D_{(r, t)}^\text{(QSS)} \equiv 
\Big\{(r, t) \Big| r > r_0, \; (-r_0 + c t ) >r \Big\} 
\end{dcases} 
\label{self_mask_c4}
\end{equation}
Outside the forcing shell, the effect of the force is zero. 

\subsection{Masking in the Exterior Using a Thin Shell of Force}
We consider a Gaussian acoustic source described by the force 
\begin{equation}
F^{(s)}(r, t) = a_s \frac{1}{d^3} \rho(\frac{r}{d}) \sin(\omega t) 
\label{source_F1}
\end{equation}
To mask this acoustic source, we introduce a thin shell of concentrated force 
surrounding the Gaussian source, given by
\begin{equation}
F^{(m)}(r, t) = a_m \delta(r-r_0) \sin(\omega t) 
\label{interf_F2}
\end{equation}
Here we adopt the convention that the acoustic source is denoted by $F^{(s)}$ 
and the interference force for masking is denoted by $F^{(m)}$. 
The quasi-steady periodic solution excited by the total force is obtained 
by superposing the solution for $F^{(s)}(r, t)$ given in \eqref{u_normal_QSS}
and that for $F^{(m)}(r, t)$ given in \eqref{u_del_r0_QSS}. 
The superposition yields 
\begin{equation}
u(r, t) = \Big(\frac{a_s}{2c^2}e^{-2\pi^2\frac{d^2}{\lambda^2}}
+ a_m \frac{r_0 \lambda}{2\pi c^2} \sin(2\pi \frac{r_0}{\lambda}) \Big)
\frac{1}{r}\sin(2\pi \frac{(c t- r)}{\lambda}) 
\label{u_F1_F2}
\end{equation}
where $\lambda \equiv \frac{2 \pi c}{\omega}$ is the wavelength.
To mask $F^{(s)}(r, t)$, we select the radius of the interference forcing shell 
$r_0 = \lambda \frac{2k+1}{4}$ where $k = $ positive integer. 
For this radius, we have
\[ r_0 = \lambda \frac{2k+1}{4}, \qquad 2\pi \frac{r_0}{\lambda} = (k+\frac{1}{2})\pi, 
\qquad \sin(2\pi \frac{r_0}{\lambda}) = (-1)^k \]
The amplitude $a_m$ of the masking force is determined by setting 
the oscillation amplitude to zero in \eqref{u_F1_F2}. The equation for $a_m$ is 
\[ \frac{a_s}{2c^2}e^{-2\pi^2\frac{d^2}{\lambda^2}}
+ a_m \frac{r_0 \lambda}{2\pi c^2} \sin(2\pi \frac{r_0}{\lambda}) = 0 \]
With the radius $r_0$ selected above, the masking amplitude $a_m$ and the wave solution are 
\begin{equation}
\begin{dcases}
a_m = (-1)^{k+1} \frac{a_s \pi}{r_0 \lambda} e^{-2\pi^2\frac{d^2}{\lambda^2}} \\[1ex] 
u(r, t) \equiv 0 \quad 
\text{for all } (r, t) \in D_{(r, t)}^\text{(QSS)} \equiv 
\Big\{(r, t) \Big| r > r_0, \; (-r_0 + c t ) >r \Big\} 
\end{dcases} \label{mask_c4}
\end{equation}
Outside the effective spatial core of the Gaussian forcing, the masking force 
effectively cancels the contribution of the source, leading to 
a vanishing net wave excitation in the quasi-steady regime. 
The radius of the masking force shell, $r_0$, is selected from the sequence of $\frac{1}{4}\lambda $, 
 $\frac{3}{4}\lambda $, $\frac{5}{4}\lambda $, $\frac{7}{4}\lambda, \ldots $. 
For a 100 Hz vibration in water, the wavelength is about $\lambda = 15\, \text{m}$.
In the expression of $F^{(m)}(r, t)$, the amplitude $a_m$ may be negative, which is 
implemented with a phase shift of $\pi$ in the time profile:
\[ -a_m \sin(\omega t) = a_m \sin(\omega t + \pi) \]

\subsection{Masking in a Small Neighborhood Using One Point-Force 
} \label{one_point_mask}
Consider the Gaussian acoustic source $F^{(s)}(|\vec{x}|, t)$ given in \eqref{source_F1}. 
In the previous subsection, we employed a thin spherical shell of force to mask  
the effect of the Gaussian source in the entire domain outside its  
spatial forcing core. 
From the perspective of operation, a full spherical shell of force 
surrounding the acoustic source is difficult to implement in practice.
In this subsection, we instead consider using a single point-force located at $\vec{x}_m$ to mask 
the effect of the Gaussian source in a small neighborhood around a specified position
$\vec{x}_d$, where possible detection sensors are located. The overall goal is to mask 
the effect of the acoustic source from detection. 

There are two key differences between the problem setup of this subsection and 
those of previous subsections: i) the objective is to mask the effect of the acoustic source 
only in a given small neighborhood, not in the infinite domain, and 
ii) we use only a single point-force for masking. 
We examine an interference point-force of the form  
\begin{equation}
F^{(m)}(\vec{x}, t) = a_m \delta(\vec{x}-\vec{x}_m) \sin(\omega t + \varphi_m) 
\label{interf_F3}
\end{equation}
where $a_m$ is the amplitude, and $\varphi_m$ is a phase shift introduced to optimize destructive interference in a small neighborhood around the detection point.

The quasi-steady periodic solution excited by the total force is obtained 
by superposing the solution for $F^{(s)}(r, t)$ given in \eqref{u_normal_QSS} 
and that for $F^{(m)}(r, t)$ given in \eqref{u_delta_QSS}. The result is
\begin{equation}
\begin{dcases}
u(\vec{x}, t) = u^{(s)}(\vec{x}, t) + u^{(m)}(\vec{x}, t) \\
\quad = \tilde{a}_s \frac{1}{|\vec{x}|}\sin(2\pi \frac{(c t- |\vec{x}|)}{\lambda})
+ \tilde{a}_m \frac{1}{|\vec{x} - \vec{x}_m|}
\sin(2\pi \frac{(c t- |\vec{x} - \vec{x}_m|)}{\lambda}+\varphi_m) \\[1ex]
\qquad \tilde{a}_s \equiv \frac{a_s}{2c^2}e^{-2\pi^2\frac{d^2}{\lambda^2}}, 
\qquad \tilde{a}_m \equiv \frac{a_m}{4 \pi c^2} 
\end{dcases} \label{u_F1_F3}
\end{equation}
where $\tilde{a}_s$ and $\tilde{a}_m$ defined above are scaled amplitudes. 
To achieve the goal of masking in a small neighborhood around $\vec{x}_d$, 
we impose the condition $ u(\vec{x}, t) \big|_{\vec{x}=\vec{x}_d} = 0$ for all $t$, 
which leads to the following requirements on the phase shift $ \varphi_m$ and 
on the scaled amplitude $\tilde{a}_m$ of the masking force:
\begin{equation}
 \varphi_m(\vec{x}_m) = \pi+ 2\pi \frac{|\vec{x}_d - \vec{x}_m|-|\vec{x}_d| }{\lambda} 
, \qquad
\tilde{a}_m(\vec{x}_m) = \tilde{a}_s \frac{|\vec{x}_d - \vec{x}_m|}{|\vec{x}_d|} 
\label{a3_phi0}
\end{equation}
Here we extended the notation to include explicitly the dependence of 
$(\varphi_m, \tilde{a}_m)$ on the location of the masking force $\vec{x}_m $.
In a small neighborhood of $\vec{x} = \vec{x}_d$, linear approximation gives 
\begin{equation}
u(\vec{x}, t) = \underbrace{u(\vec{x}, t) \big|_{\vec{x}=\vec{x}_d}}_{=0}
+\Big\langle \nabla_{\vec{x}} u(\vec{x}, t) \big|_{\vec{x}=\vec{x}_d}, \; 
(\vec{x}-\vec{x}_d) \Big\rangle + \cdots \label{linear_eq}
\end{equation}
It is clear in \eqref{a3_phi0} that the condition $ u(\vec{x}, t) \big|_{\vec{x}=\vec{x}_d} = 0$ 
can be satisfied by placing a point-force at any location $\vec{x}_m$, provided that its phase shift 
$\varphi_m$ and scaled amplitude $\tilde{a}_m$ are chosen according to \eqref{a3_phi0}. 
To minimize $u(\vec{x}, t)$ in a small neighborhood of $\vec{x}_d$, we minimize the gradient 
$\nabla_{\vec{x}} u(\vec{x}, t) \big|_{\vec{x}=\vec{x}_d}$, thereby reducing the leading-order variation of the field around $\vec{x}_d$. 
This ensures that not only is the field canceled exactly at $\vec{x}_d$, but that nearby points 
also experience minimal residual wave amplitude due to reduced spatial variation. 

To calculate $\nabla_{\vec{x}} u(\vec{x}, t) \big|_{\vec{x}=\vec{x}_d}$, 
we differentiate the two terms in \eqref{u_F1_F3}, and substitute in 
$\varphi_m(\vec{x}_m) $ and $\tilde{a}_m(\vec{x}_m)$ from \eqref{a3_phi0}. 
The gradient of the Gaussian source term is 
\begin{align*}
& \nabla_{\vec{x}}\Big( \frac{1}{|\vec{x}|}
\sin(2\pi \frac{(c t- |\vec{x}|)}{\lambda}) \Big) \Big|_{\vec{x} = \vec{x}_d} \\
& \qquad = \frac{-\vec{x}_d}{|\vec{x}_d|^3}\sin(2\pi \frac{(c t- |\vec{x}_d|)}{\lambda}) 
- \frac{2 \pi \vec{x}_d}{\lambda |\vec{x}_d|^2}\cos(2\pi \frac{(c t- |\vec{x}_d|)}{\lambda})
\end{align*}
Next, the gradient of the masking term is 
\begin{align*}
& \nabla_{\vec{x}}\Big( \frac{1}{|\vec{x} - \vec{x}_m|}
\sin(2\pi \frac{(c t- |\vec{x} - \vec{x}_m|)}{\lambda}+\varphi_m)
\Big) \Big|_{\vec{x} = \vec{x}_d} \\
& \qquad = \frac{-(\vec{x}_d-\vec{x}_m)}{|\vec{x}_d-\vec{x}_m|^3}
\sin(2\pi \frac{(c t- |\vec{x}_d - \vec{x}_m|)}{\lambda}+\varphi_m) 
\nonumber \\
& \qquad \qquad -\frac{2 \pi (\vec{x}_d-\vec{x}_m)}{\lambda |\vec{x}_d-\vec{x}_m|^2}
\cos(2\pi \frac{(c t- |\vec{x}_d - \vec{x}_m|)}{\lambda}+\varphi_m)  \\
& \qquad = \frac{\vec{x}_d-\vec{x}_m}{|\vec{x}_d-\vec{x}_m|^3}
\sin(2\pi \frac{(c t- |\vec{x}_d|)}{\lambda}) 
+\frac{2 \pi (\vec{x}_d-\vec{x}_m)}{\lambda |\vec{x}_d-\vec{x}_m|^2}
\cos(2\pi \frac{(c t- |\vec{x}_d|)}{\lambda}) 
\end{align*}
Using the two gradient expressions above, we obtain $\nabla_{\vec{x}} u(\vec{x}, t)$ 
at $\vec{x}_d$:
\begin{align}
& \nabla_{\vec{x}} u(\vec{x}, t) \big|_{\vec{x}=\vec{x}_d} = 
 \frac{\tilde{a}_s}{\lambda |\vec{x}_d|}
\Big(\frac{\lambda}{|\vec{x}_d-\vec{x}_m|}\widehat{(\vec{x}_d-\vec{x}_m)}
-\frac{\lambda}{|\vec{x}_d|}\widehat{\vec{x}_d} \Big) 
\sin(2\pi \frac{(c t- |\vec{x}_d|)}{\lambda}) \nonumber \\[1ex]
& \qquad\qquad + \underbrace{\frac{2 \pi \tilde{a}_s}{\lambda |\vec{x}_d|} 
\Big(\widehat{(\vec{x}_d-\vec{x}_m)}-\widehat{\vec{x}_d} \Big) 
\cos(2\pi \frac{(c t- |\vec{x}_d|)}{\lambda})}_{\text{dominant part}}
\label{Delta_u_1}
\end{align}
Here, $\widehat{\vec{w}}\equiv \frac{\vec{w}}{|\vec{w}|}$ denotes the unit vector 
in the direction of $\vec{w}$. 
We now consider the regime of $\lambda / |\vec{x}_d| \ll 1$ in which the distance between
the Gaussian source and the detection sensor is much larger than the wavelength. 
In this case, the dominant part of $\nabla_{\vec{x}} u(\vec{x}, t) \big|_{\vec{x}}$
is the second term, as marked in \eqref{Delta_u_1}. 
To minimize the residual field near $\vec{x}_d$,  we set the dominant part to zero: 
\[ 
\widehat{(\vec{x}_d-\vec{x}_m)}-\widehat{\vec{x}_d} = 0
\]
which implies that the acoustic source, the masking force and the detection sensor are aligned on one line in $\mathbb{R}^3$. 
\begin{equation}
\vec{x}_m = \beta \vec{x}_d, \;\; \text{with} \;\; \beta < 1
\label{x0_cond_1}
\end{equation}
In other words, the masking point-force should be placed along the radial line connecting 
the acoustic source and the detection sensor.

Under geometric constraint (\ref{x0_cond_1}), the remaining term in the gradient is 
\begin{align}
& \nabla_{\vec{x}} u(\vec{x}, t) \big|_{\vec{x}=\vec{x}_d} = 
\tilde{a}_s \frac{\widehat{\vec{x}_d}}{ |\vec{x}_d|^2} \big(\frac{\beta}{1-\beta}\big) 
\sin(2\pi \frac{(c t- |\vec{x}_d|)}{\lambda}) 
\label{Delta_u_rem}
\end{align}
Even with constraint (\ref{x0_cond_1}), parameter $\beta$ remains free. 
To further suppress the residual field near $\vec{x}_d$, we minimize 
 $\nabla_{\vec{x}} u(\vec{x}, t) \big|_{\vec{x}=\vec{x}_d}$ in \eqref{Delta_u_rem} 
by selecting $\beta $ as small as possible (i.e., placing the masking point-force
as close as possible to the Gaussian source). 
Mathematically, if the masking force is placed right at the center of the Gaussian source 
($\beta = 0$), it will completely cancel the Gaussian source outside 
its spatial forcing core. 
Operationally, the region near the center of the Gaussian source may be inaccessible. 
Suppose the ball $ B(0, \varepsilon_s)$ around the Gaussian source is excluded when 
placing the masking force. Then, the location of the masking force is restricted by 
$|\vec{x}_m| \ge \varepsilon_s$ and parameter $\beta $ in $\vec{x}_m = \beta \vec{x}_d $
is restricted by $|\beta| \ge \varepsilon_s/ |\vec{x}_d|$. 
Combining \eqref{interf_F3}, \eqref{a3_phi0} and \eqref{x0_cond_1}, 
we obtain the optimal point-force for masking the Gaussian source in a small neighborhood of $\vec{x}_d$. 
\begin{equation}
\begin{dcases}
F^{(m)}(\vec{x}, t) = a_m \delta(\vec{x}-\vec{x}_m) \sin(\omega t + \varphi_m) 
\\[1ex]
\qquad \vec{x}_m = \beta \vec{x}_d = \pm \varepsilon_s \hat{\vec{x}}_d, 
\quad \beta = \pm \varepsilon_s/|\vec{x}_d|, 
\\[1ex]
\qquad \varphi_m = \pi (1-\frac{2\beta |\vec{x}_d| }{\lambda}) 
 = \pi (1 \mp \frac{2\varepsilon_s}{\lambda}) 
 \\[1ex] 
\qquad a_m = a_s (1-\beta) 2\pi e^{-2\pi^2\frac{d^2}{\lambda^2}} 
\end{dcases}
\label{optimal_pt_frc}
\end{equation}
Figure \ref{fig_01} illustrates the optimal location of a point-force for masking 
the effect of the Gaussian source in a small region of possible detection sensors. 
This optimal location shifts in response to changes in the prescribed region of sensors.
\begin{figure}[!h]
\vskip 0.5cm
\begin{center}
\psfig{figure=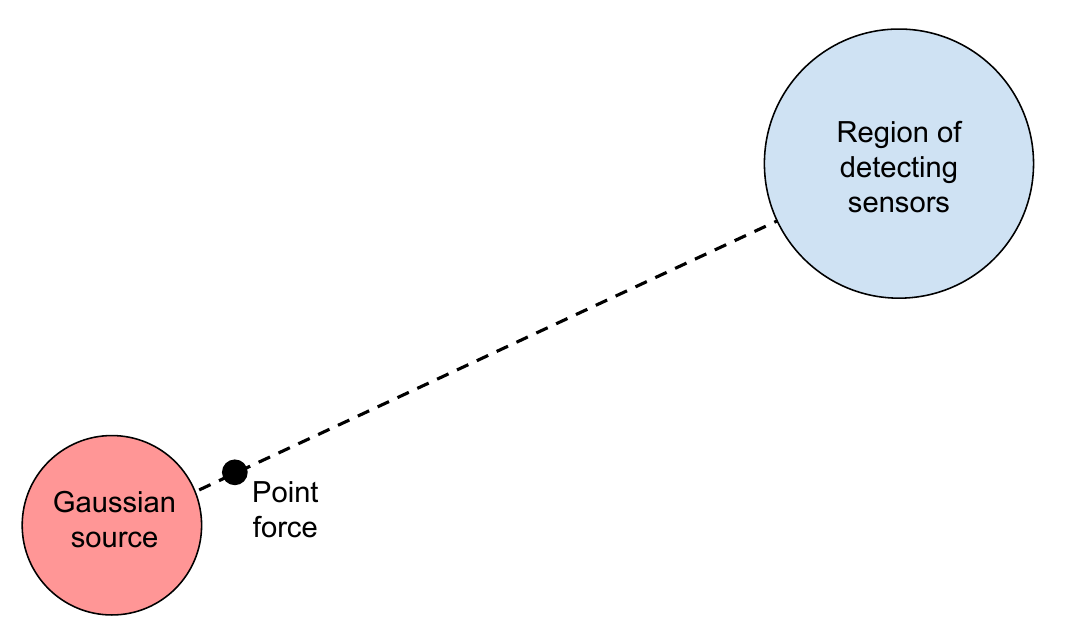, height=2.5in}
\end{center}
\vskip -0.75cm
\caption{Optimal location of a point-force for masking the effect of a known 
Gaussian source in a small region around a given location $\vec{x}_d$.}
\label{fig_02A}
\end{figure}

We investigate the efficacy of the optimal point-force \eqref{optimal_pt_frc}
in masking the Gaussian source.
Let $B(\vec{x}_d, \varepsilon_d)$ denote the ball of radius $\varepsilon_d$ 
centered at point $\vec{x}_d$, representing the small region of sensors. 
We analyze the signal within this region both before and after the application of 
the optimal point-force in masking.
Let $\text{amp}\{u(\vec{x}, \cdot)\}$ denote the oscillation amplitude of 
$u(\vec{x}, t)$ over time $t$ at location $\vec{x}$. 
\[ \text{amp}\{u(\vec{x}, \cdot)\} \equiv \max_{t} |u(\vec{x}, t)| \]
Prior to applying the optimal point-force, the signal amplitude is approximately 
\begin{equation}
 \text{amp}\{u^{(s)}(\vec{x}, \cdot)\}
\approx \tilde{a}_s \frac{1}{|\vec{x}_d|}, \quad 
\vec{x} \in B(\vec{x}_d, \varepsilon_d) 
\label{amp_u_signal}
\end{equation}
Upon applying the optimal point-force \eqref{optimal_pt_frc}, and utilizing the expressions (\ref{linear_eq}) and (\ref{Delta_u_rem}), we observe that
the signal amplitude becomes
\begin{align}
& \text{amp}\{u(\vec{x}, \cdot)\}
 \approx \Big| \tilde{a}_s \frac{\langle \widehat{\vec{x}_d}, 
(\vec{x}-\vec{x}_d) \rangle}{ |\vec{x}_d|^2} \big(\frac{\beta}{1-\beta}\big) 
\sin(2\pi \frac{(c t- |\vec{x}_d|)}{\lambda}) \Big|, 
\quad \vec{x} \in B(\vec{x}_d, \varepsilon_d) \nonumber \\[1ex]
& \hspace*{2cm} \le \tilde{a}_s \frac{\varepsilon_d}
{ |\vec{x}_d|^2} \big(\frac{|\beta| }{1-\beta}\big), 
\quad \vec{x} \in B(\vec{x}_d, \varepsilon_d) 
\label{amp_u_mask}
\end{align} 
At each position $\vec{x}$, we define the normalized residual amplitude as the ratio of the signal amplitudes after and before applying the optimal point-force:  
\begin{equation}
A^\text{(NR)}(\vec{x}) \equiv 
\frac{\underbrace{\;\;\text{amp}\{u(\vec{x}, \cdot)\} \;\;}_\text{after masking}}{\;\;\underbrace{\;\;\text{amp}\{u^{(s)}(\vec{x}, \cdot)\}\;\;}_\text{ before masking}\;\;} 
\label{resi_amp_def}
\end{equation}
We use this normalized residual amplitude as a metric in assessing the effectiveness of masking.  
Within the region $B(\vec{x}_d, \varepsilon_d)$, employing the expressions \eqref{amp_u_signal} 
and \eqref{amp_u_mask} for the two amplitudes, we obtain the estimate: 
\begin{equation}
\boxed{\quad 
A^\text{(NR)}(\vec{x})
= \frac{\text{amp}\{u(\vec{x}, \cdot)\}}
{\text{amp}\{u^{(s)}(\vec{x}, \cdot)\}} 
\le \frac{\varepsilon_d}{ |\vec{x}_d|} \big(\frac{|\beta| }{1-\beta}\big), 
\quad \vec{x} \in B(\vec{x}_d, \varepsilon_d) }
\label{resi_amp_2}
\end{equation}
In summary, the optimal point-force given in \eqref{optimal_pt_frc} is engineered 
to mask the effect of the Gaussian source in a small neighborhood of the given location 
$\vec{x}_d$. The optimal point-force varies with the prescribed location 
$\vec{x}_d$. 

\subsection{Masking in a Small Neighborhood Using Two Point-Forces}
\label{two_point_mask}
Consider the Gaussian acoustic source $F^{(s)}(|\vec{x}|, t)$ described in \eqref{source_F1}.
In the preceding subsection, we derived the optimal solution when a single point-force
is used to mask the effect of $F^{(s)}(|\vec{x}|, t)$ in a small neighborhood of $\vec{x}_d$.
The approach is to first make the dominant part in gradient 
$\nabla_{\vec{x}} u(\vec{x}, t) \big|_{\vec{x}=\vec{x}_d}$ zero and then 
minimize the remaining part. With a single point-force, we cannot achieve 
$\nabla_{\vec{x}} u(\vec{x}, t) \big|_{\vec{x}=\vec{x}_d} = 0$ for all $t$.

In this subsection, we use two point-forces to enhance the masking effect. 
With two point-forces, we are able to make 
$\nabla_{\vec{x}} u(\vec{x}, t) \big|_{\vec{x}=\vec{x}_d} = 0$ for all $t$.

First we apply each of the two point-forces separately in two different maskings, each 
utilizing only a single point-force. In these two separate maskings, we place each of 
the two point-forces according to the optimal solution obtained in \eqref{optimal_pt_frc}
for a single point-force masking.
The first point-force is assigned $\beta_1 > 0$ while the second point-force is assigned 
$\beta_2 < 0$.
\begin{equation}
\begin{dcases}
\text{Force 1: } F^{(m,1)}(\vec{x}, t) = a_m^{(1)} \delta(\vec{x}-\vec{x}_m^{(1)}) 
\sin(\omega t + \varphi_m^{(1)}) 
\\[1ex]
\qquad \vec{x}_m^{(1)} = \beta_1 \vec{x}_d = \varepsilon_s \widehat{\vec{x}_d}, 
\qquad \beta_1 = \varepsilon_s/|\vec{x}_d| \equiv \beta_s, 
\\[1ex]
\qquad \varphi_m^{(1)} = \pi (1-\frac{2\beta_1 |\vec{x}_d| }{\lambda}) 
=  \pi (1-\frac{2\varepsilon_s }{\lambda})
 \\[1ex] 
\qquad a_m^{(1)} = a_s (1-\beta_1) 2\pi e^{-2\pi^2\frac{d^2}{\lambda^2}} 
\end{dcases}
\label{optimal_pt_frc_1}
\end{equation}
\begin{equation}
\begin{dcases}
\text{Force 2: } F^{(m,2)}(\vec{x}, t) = a_m^{(2)} \delta(\vec{x}-\vec{x}_m^{(2)}) 
\sin(\omega t + \varphi_m^{(2)}) 
\\[1ex]
\qquad \vec{x}_m^{(2)} = \beta_2 \vec{x}_d = -\varepsilon_s \widehat{\vec{x}_d},
\qquad \beta_2 = -\varepsilon_s/|\vec{x}_d| = -\beta_s, 
\\[1ex]
 \qquad \varphi_m^{(2)} = \pi (1-\frac{2\beta_2 |\vec{x}_d| }{\lambda}) 
 = \pi (1+\frac{2\varepsilon_s }{\lambda})
 \\[1ex] 
\qquad a_m^{(2)} = a_s (1-\beta_2) 2\pi e^{-2\pi^2\frac{d^2}{\lambda^2}} 
\end{dcases}
\label{optimal_pt_frc_2}
\end{equation}
It is worthwhile to emphasize that the two point-forces described in \eqref{optimal_pt_frc_1} 
and \eqref{optimal_pt_frc_2} are the optimal solutions in two separate maskings, 
each utilizing only a single point-force. 
In each separate masking, the resulting solution $u^{(j)}(\vec{x}, t)$ at the sensor location 
$\vec{x}_d$ is zero for all $t$: 
$ u^{(j)}(\vec{x}, t) \big|_{\vec{x}=\vec{x}_d} = 0, \; j = 1, 2 $. 
In the resulting gradient $\nabla_{\vec{x}} u^{(j)}(\vec{x}, t) \big|_{\vec{x}=\vec{x}_d}, \; j = 1, 2$
of each separate masking, the dominant $\cos(\;)$ part vanishes for all $t$, guaranteed by 
the location of each point-force, $\vec{x}_m^{(j)}, \; j = 1, 2$, in \eqref{optimal_pt_frc_1} 
and \eqref{optimal_pt_frc_2}. 
The remaining gradient of each separate masking, as given 
in \eqref{Delta_u_rem}, is proportional to 
\[ \nabla_{\vec{x}} u^{(j)}(\vec{x}, t) \big|_{\vec{x}=\vec{x}_d} \propto 
(\frac{\beta_j}{1-\beta_j}) \sin(2\pi \frac{(c t- |\vec{x}_d|)}{\lambda}), \quad j = 1, 2. \]
Now we use a weighted average to combine the two optimal solutions 
from two separate maskings to construct a new masking configuration that utilizes both of the two forces. 
In the weighted average, we use $\gamma$ fraction of force 1 and 
$(1-\gamma)$ fraction of force 2.
\begin{equation}
F^{(m)}(\vec{x}, t) = \gamma F^{(m,1)}(\vec{x}, t)+(1-\gamma) F^{(m,2)}(\vec{x}, t)
\label{F_wt_avg_1}
\end{equation}
We apply the weighted average $F^{(m)}(\vec{x}, t)$ to mask the Gaussian source. 
The residual gradient $\nabla_{\vec{x}} u(\vec{x}, t) \big|_{\vec{x}=\vec{x}_d} 
= \nabla_{\vec{x}} u^{(1)}(\vec{x}, t) \big|_{\vec{x}=\vec{x}_d}
+\nabla_{\vec{x}} u^{(2)}(\vec{x}, t) \big|_{\vec{x}=\vec{x}_d} $ 
is proportional to
\[ \nabla_{\vec{x}} u(\vec{x}, t) \big|_{\vec{x}=\vec{x}_d} \propto 
\Big(\gamma \frac{\beta_1}{1-\beta_1}+(1-\gamma) \frac{\beta_2}{1-\beta_2}\Big) 
\sin(2\pi \frac{(c t- |\vec{x}_d|)}{\lambda}) \]
To make the residual gradient exactly zero at $\vec{x}_d$, we set 
$\gamma \frac{\beta_1}{1-\beta_1}+(1-\gamma) \frac{\beta_2}{1-\beta_2}=0 $.
Using $\beta_1 = \beta_s $ and $\beta_2 = -\beta_s$ from \eqref{optimal_pt_frc_1} 
and \eqref{optimal_pt_frc_2}, we obtain an equation on $\gamma$. 
\begin{equation}
\gamma \frac{\beta_s}{1-\beta_s}+(1-\gamma) \frac{(-\beta_s)}{1+\beta_s}=0
\label{b1_b2_eq}
\end{equation}
Solving for coefficient $\gamma$ in \eqref{b1_b2_eq} yields $\gamma = (1-\beta_s)/2 $. 
With this value of $\gamma$, the weighted average $F^{(m)}(\vec{x}, t)$ 
makes the residual gradient exact zero at $\vec{x}_d$, and 
based on the Taylor linear approximation, it minimizes the residual solution 
in a small neighborhood of $\vec{x}_d$. 
Therefore, the optimal weighted average of the two point-forces for masking 
the effect of the Gaussian source in a small neighborhood of $\vec{x}_d$ is
\begin{equation}
\boxed{\quad F^{(m)}(\vec{x}, t) = \frac{1-\beta_s}{2} F^{(m,1)}(\vec{x}, t)
+\frac{1+\beta_s}{2}F^{(m,2)}(\vec{x}, t), \quad 
\beta_s \equiv \frac{\varepsilon_s}{|\vec{x}_d|}  \quad} 
\label{F_wt_avg_2}
\end{equation}
where $F^{(m,1)}(\vec{x}, t)$ and $F^{(m,2)}(\vec{x}, t)$ are optimal solutions
given in \eqref{optimal_pt_frc_1} and \eqref{optimal_pt_frc_2}, respectively,
in two separate maskings, each utilizing only one point-force. 

In summary, for masking the effect of a Gaussian source in a given small neighborhood 
using two point-forces, the optimal design is obtained through the following steps:
\begin{enumerate}
\item Individual optimization: optimize each point-force in a separate masking that uses only one point-force:
 \begin{itemize}
 \item Force 1 is optimized for $\beta_1 > 0$, as detailed in \eqref{optimal_pt_frc_1}.
 \item Force 2 is optimized for $\beta_2 < 0$, as specified in 
 \eqref{optimal_pt_frc_2}.
 \end{itemize}
\item Weighted average combination: Use a weighted average of the two individually optimized 
forces in a new masking configuration. Adjust the coefficient to make the residual 
gradient exactly zero at the given location, as described in \eqref{F_wt_avg_2}.
\end{enumerate}

\begin{figure}[!h]
\vskip 0.5cm
\begin{center}
\psfig{figure=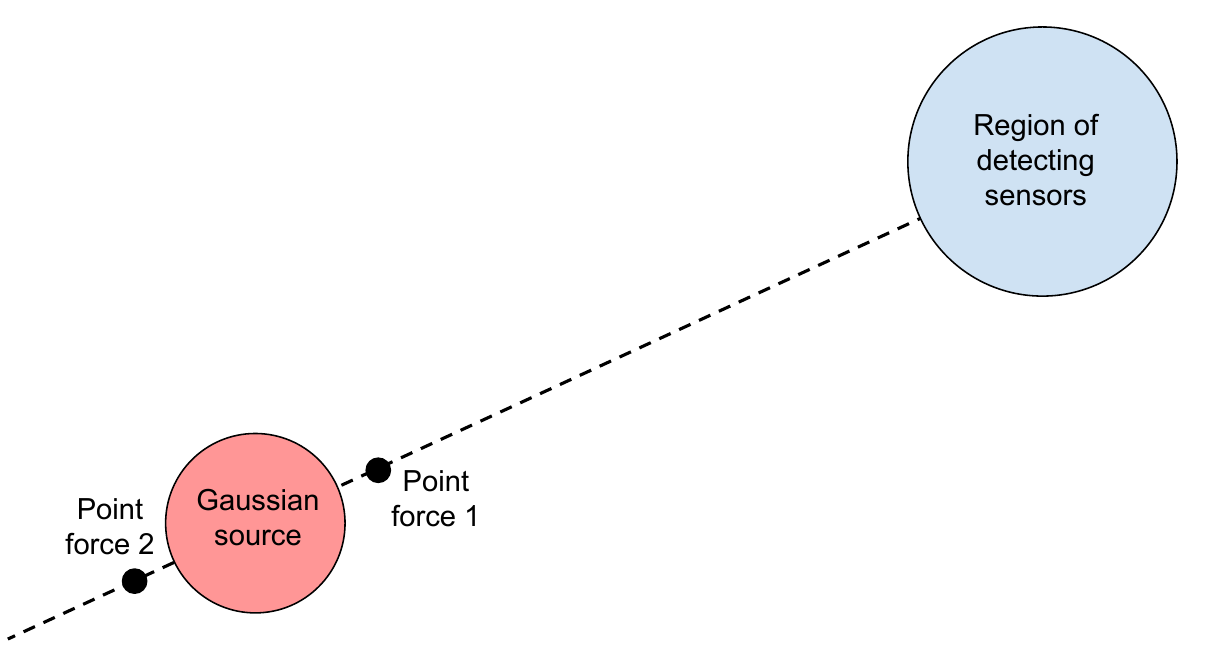, height=2.5in}
\end{center}
\vskip -0.75cm
\caption{Optimal placement of two point-forces to mask the effect of a known 
Gaussian source in a small region around a given location $\vec{x}_d$.}
\label{fig_02B}
\end{figure}
%

\subsection{Performance of Maskings Using One or Two Point-Forces} 
\begin{figure}[!h]
\vskip 0.5cm
\begin{center}
\psfig{figure=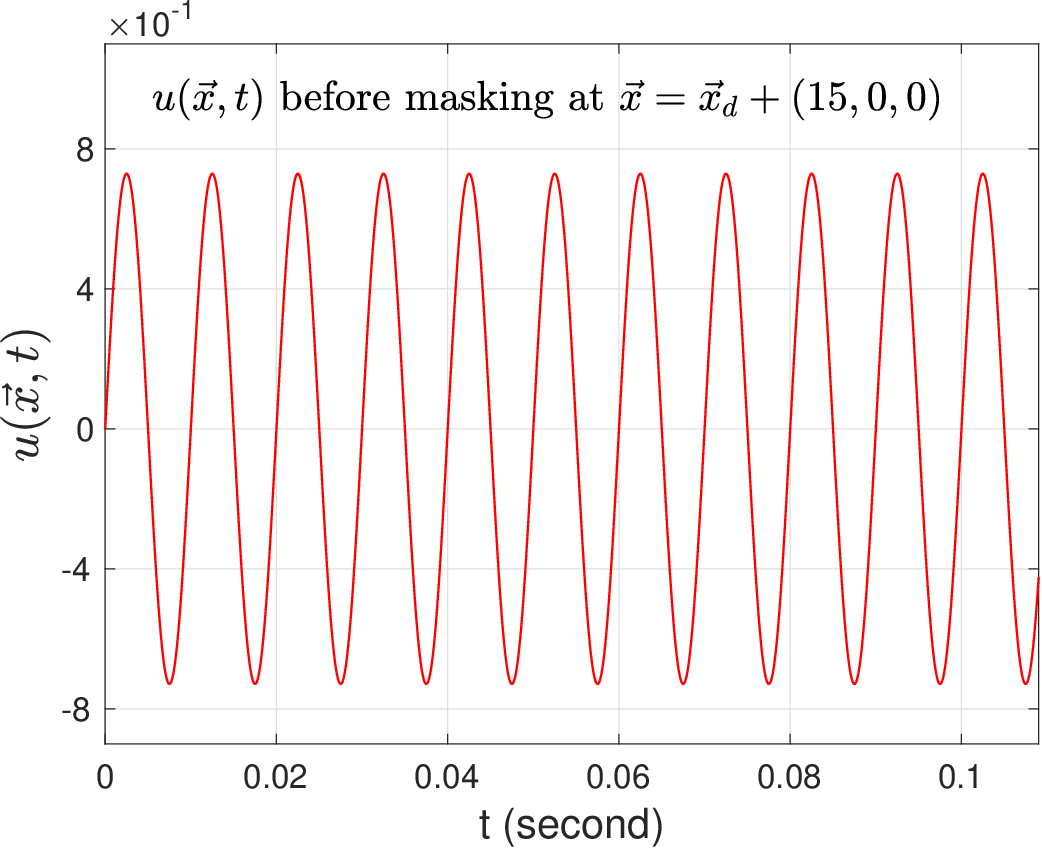, width=2.8in}\qquad
\psfig{figure=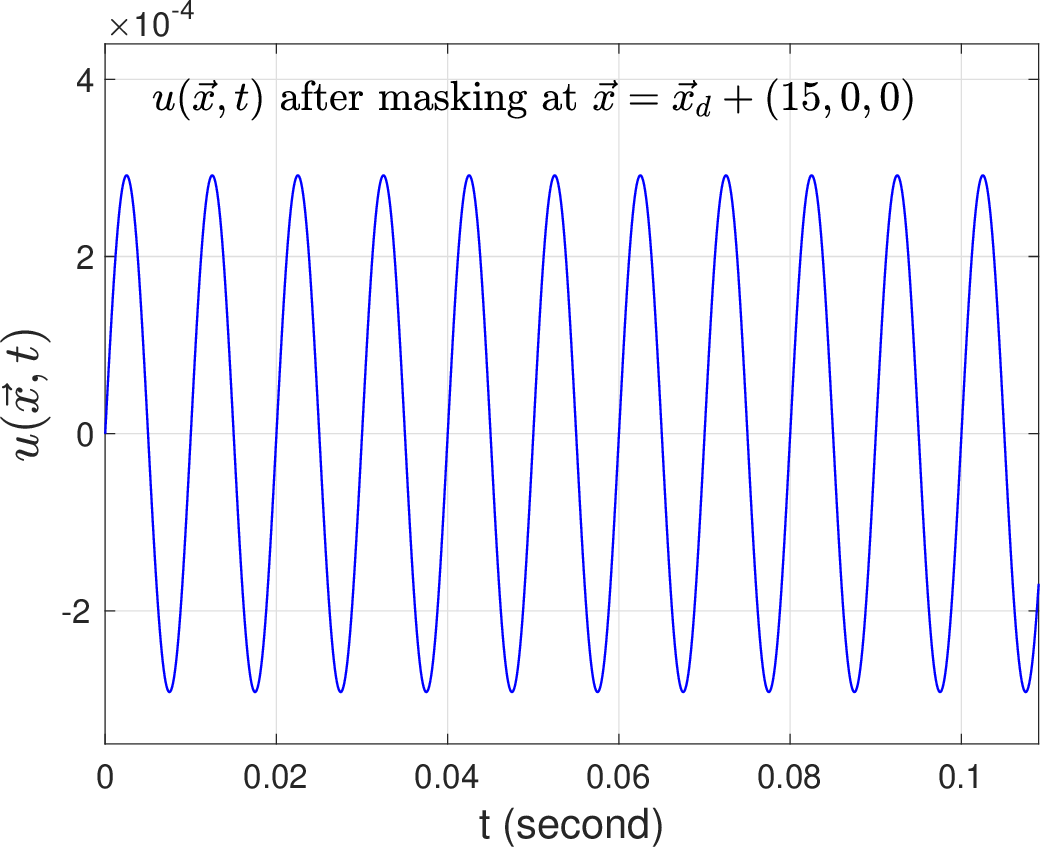, width=2.8in} \\[1ex]
\psfig{figure=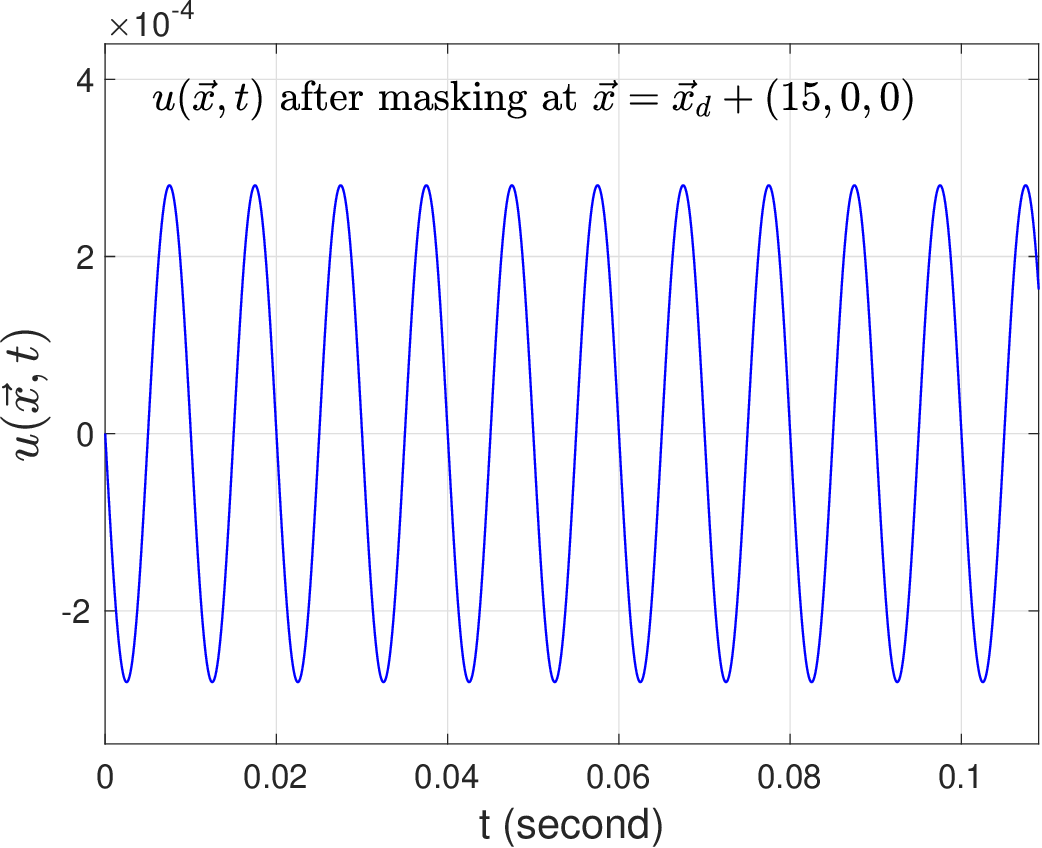, width=2.8in}\qquad
\psfig{figure=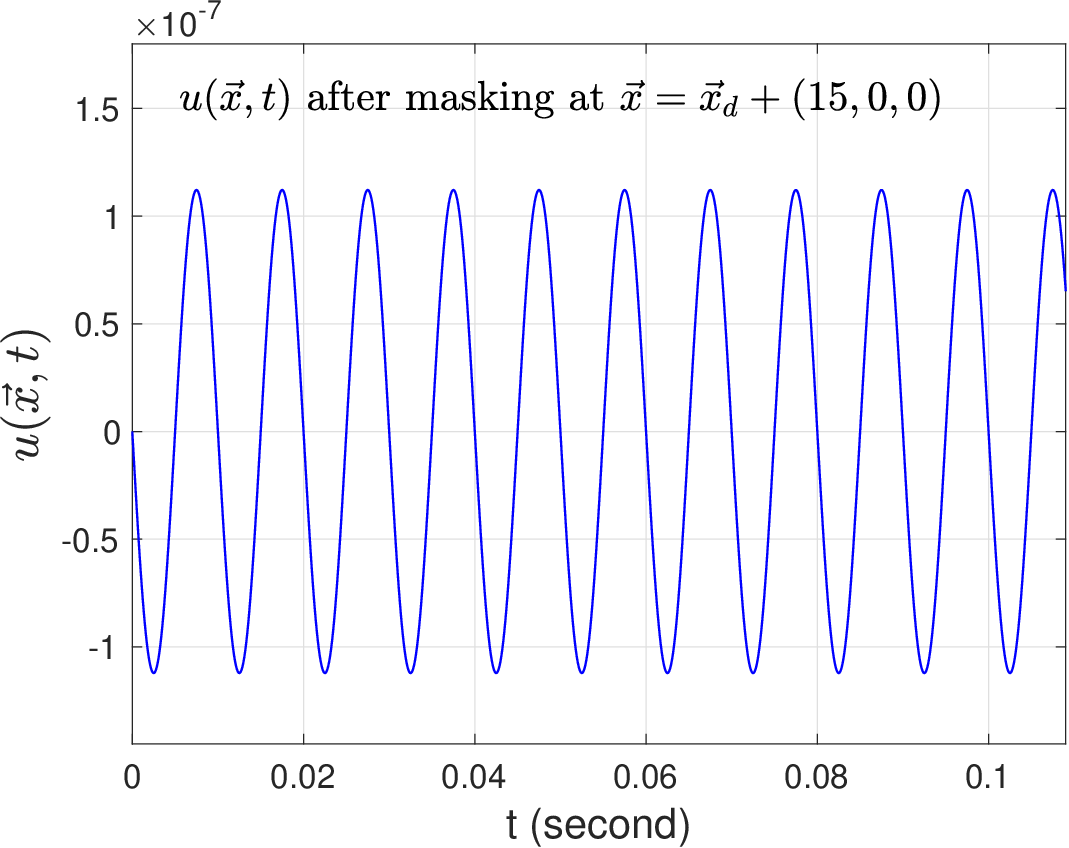, width=2.8in}
\end{center}
\vskip -0.75cm
\caption{$u(\vec{x}, t)$ vs $t$ at a fixed point before and after masking. 
Top left: source only. Top right: source + point-force 1.
Bottom left: source + point-force 2. Bottom right: source + optimal weighted average
of the two forces.}
\label{fig_03}
\end{figure}

We evaluate the effectiveness of the optimal masking solutions derived in subsections 
\ref{one_point_mask} and \ref{two_point_mask}. 
In our analysis, we use the parameters and coordinate system listed below. 
\begin{itemize}
\item Sound speed: $c = 1500\, \text{m/s}$, the speed of acoustic wave in water. 
\item Center of the Gaussian source: at the origin, $(0,0,0)$. 
\item Region where possible sensors are located: we model the region as a ball, 
$ B(\vec{x}_d, \varepsilon_d)$; 
the center is $\vec{x}_d = |\vec{x}_d|  (1, 0, 0)$ with $ |\vec{x}_d| = 750\,\text{m}$; 
the radius is $\varepsilon_d = 15\,\text{m}$. 
This is the target region in which we want to mask 
the effect of the Gaussian source from detection. 
\item Frequency and wavelength: we consider an acoustic source of $100$ Hz; 
the angular frequency is $\omega = 100{\times} 2\pi/\text{s}$; 
the wavelength is $\lambda = \frac{2\pi c}{\omega} = 15\,\text{m}$.
\item Inaccessible region around the Gaussian source: we model the region as a ball, 
$ B(0, \varepsilon_s)$; the center is the origin; the radius is $\varepsilon_s = 15\,\text{m}$. 
This is the excluded region for positioning point-force(s) when we design masking configurations. 
\end{itemize}
We examine the performance of three optimal masking configurations: 
each of the first two configurations uses only one point-force for interference, 
respectively with $\beta > 0$ and with $\beta < 0$, 
as detailed in \eqref{optimal_pt_frc_1} and \eqref{optimal_pt_frc_2}; 
the third configuration uses the optimal weighted average of the two point-forces for interference, 
as described in \eqref{F_wt_avg_2}.
We first examine the temporal evolution of the acoustic field $u(\vec{x}, t)$ 
at a fixed location $\vec{x} = \vec{x}_d + \varepsilon_d (1, 0, 0) $,
before and after applying each of three optimal masking configurations.
The results are shown in the four panels of Figure \ref{fig_03}.
\begin{itemize}
    \item Top left: the source field before any masking is applied.
    \item Top right: the net field when only point-force 1 is applied.
    \item Bottom left: the net field when only point-force 2 is applied.
    \item Bottom right: the net field when the weighted average of two forces is applied.
\end{itemize}
In all panels, the residual signal exhibits periodicity with
the same frequency as the forcing frequency. Notably, the amplitude of the residual signal 
is significantly reduced from that of the source signal when either point-force 1 
or point-force 2 is applied (top right and bottom left panels). 
When the optimal weighted average of the two point-forces is applied (bottom right panel), 
the residual signal amplitude is further reduced by more than three orders of magnitude (bottom right panel). 
It is important to notice that the vertical scales of the four panels are individually set to 
accommodate the amplitudes in each case. 
\begin{figure}[!h]
\vskip 0.5cm
\begin{center}
\psfig{figure=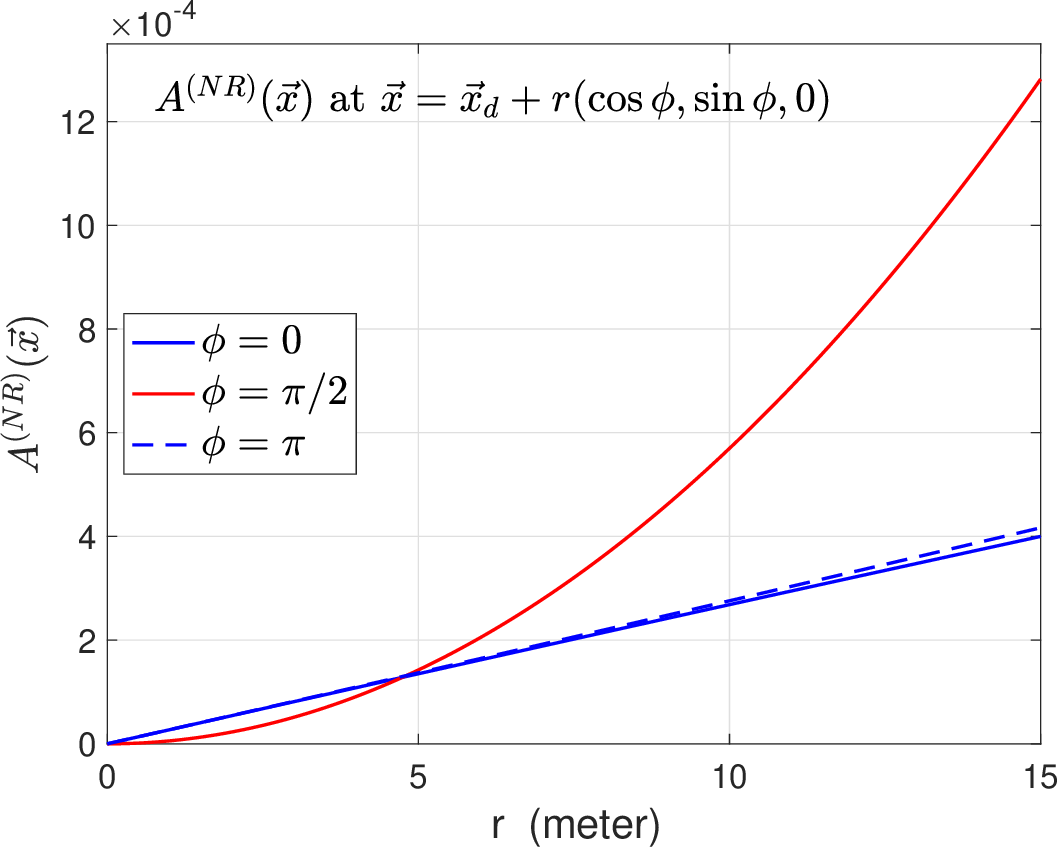, width=2.85in}\qquad
\psfig{figure=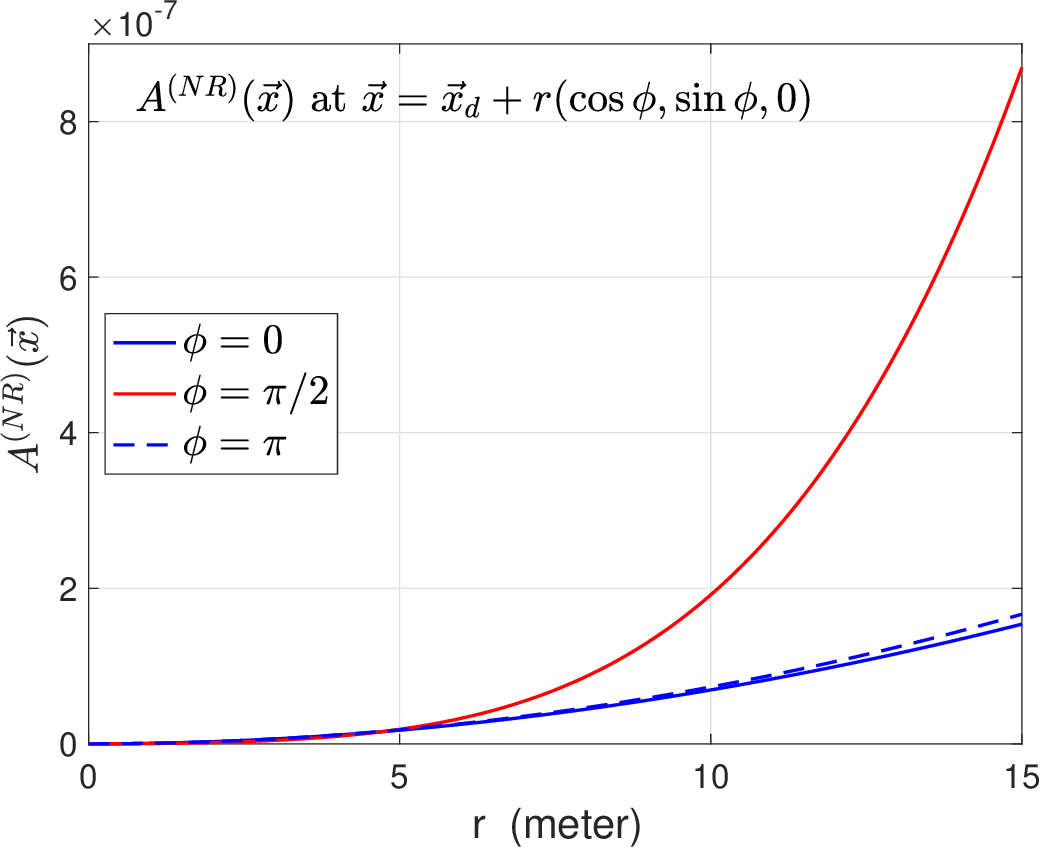, width=2.8in} 
\end{center}
\vskip -0.75cm
\caption{Normalized residual amplitude $A^\text{(NR)}(\vec{x})$
defined in \eqref{resi_amp_def}, as a function of 
radial distance $r$ along three radial lines from $\vec{x}_d$. 
Left: source + force 1. 
Right: source + the optimal weighted average of forces 1 and 2.}
\label{fig_04}
\end{figure}

Recall that, by the design of force 1 and force 2, after applying any of the three masking
configurations, the residual signal is zero for all $t$ at location $\vec{x}_d$. 
We examine the normalized residual amplitude $A^\text{(NR)}(\vec{x})$ 
defined in \eqref{resi_amp_def}, along several radial lines in the $(x,y)$-plane. 
Each radial line is originated from $\vec{x}_d$, the center of the region
where possible sensors are located, is specified by the azimuthal angle $\phi$ 
in the $(x,y)$-plane, and is parameterized by the radial distance $r$ from $\vec{x}_d$. 
We focus on the part inside the given region $B(\vec{x}_d, \varepsilon_d)$. 
\[\vec{x}(r) = \vec{x}_d + r(\cos\phi, \sin\phi, 0), \quad 0 \le r \le \varepsilon_d  \]
We consider three radial lines with azimuthal angle $\phi = 0,\, \pi/2,\, \pi$.
Figure \ref{fig_04} shows the normalized residual amplitude $A^\text{(NR)}(\vec{x})$ 
as a function of the radial distance $r$ along the three radial lines. 
We compare the results of two masking configurations: 
source + force 1 (left panel); source + the optimal weighted average of forces 1 and 2 
(right panel).

When only force 1 (or force 2) is applied in masking, the residual gradient at $\vec{x}_d$ 
is in the direction of $\vec{x}_d$ as derived in \eqref{Delta_u_rem}. 
Based on the Taylor approximation, for small $r$, the amplitude of the residual signal 
is proportional to $r$ along the $\vec{x}_d$-direction (i.e., for $\phi = 0$ and $\phi = \pi$), 
and is proportional to $r^2$ along the direction perpendicular to $\vec{x}_d$ (i.e., for $\phi=\pi/2$). 
This theoretical prediction is confirmed in the left panel of Figure \ref{fig_04}
where, for small $r$ (e.g., $r \le 4 \text{m}$), the residual amplitude along $\phi=\pi/2$ is lower than those along
$\phi = 0$ and $\phi = \pi$. However, for moderate values of $r$ (e.g., $r=10 \sim 15 \text{m}$), 
the residual amplitude along $\phi=\pi/2$ exceeds those along $\phi = 0$ and $\phi = \pi$. 
When the optimal weighted average of the two forces is applied in masking, 
the residual gradient at $\vec{x}_d$ is exactly zero, as ensured by the design 
of the optimal weights. As a result, the amplitude of the residual signal is expected to be
proportional to $r^2$ along all directions. 
In the right panel of Figure \ref{fig_04}, for small $r$, the residual amplitude scales as
$r^2$ along all three directions examined. 
However, for moderate values of $r$, along the $\phi=\pi/2$ direction, 
the $r^4$ term takes over as the dominant term, 
leading to that the residual amplitude along $\phi=\pi/2$ surpasses those along 
$\phi = 0$ and $\phi = \pi$. 
Comparing the two panels of Figure \ref{fig_04}, 
we observe that the optimal weighted average of the two forces is significantly more 
effective than either individual point-force in masking the Gaussian source. 
In a ball of radius 15m around $\vec{x}_d$, switching from a single point-force to 
the optimal weighted average of two forces in masking, 
the residual amplitude is further reduced by a factor of more than $10^3$. 

\begin{figure}[!h]
\vskip 0.5cm
\begin{center}
\psfig{figure=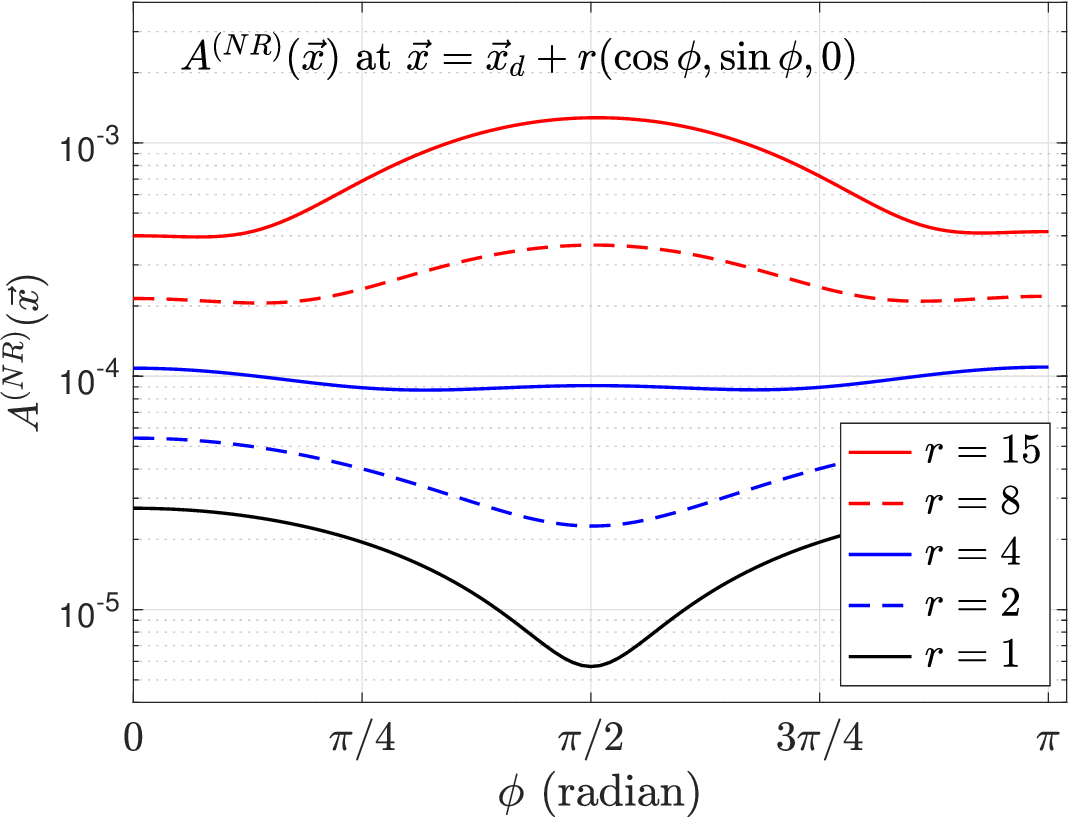, width=2.8in}\qquad
\psfig{figure=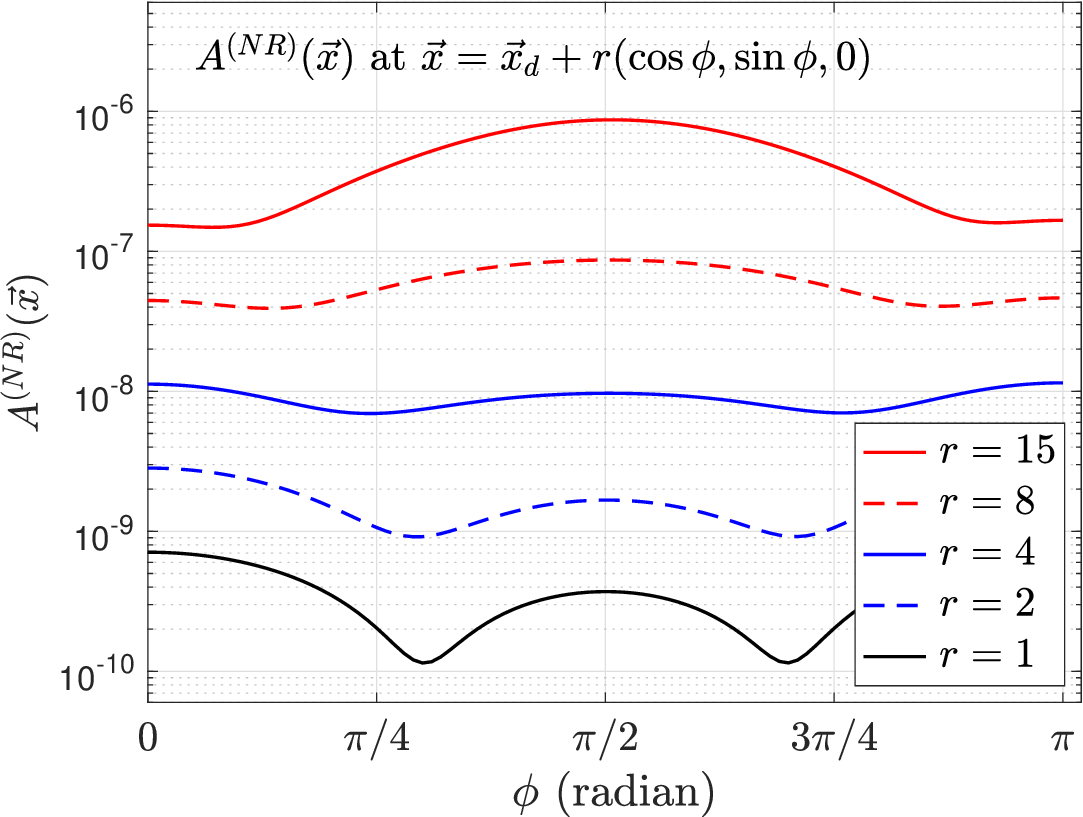, width=2.85in} 
\end{center}
\vskip -0.75cm
\caption{Normalized residual amplitude $A^\text{(NR)}(\vec{x})$
defined in \eqref{resi_amp_def}, as a function of 
azimuthal angle $\phi$ along five semi-circles around $\vec{x}_d$. 
Left: source + force 1. 
Right: source + the optimal weighted average of forces 1 and 2.}
\label{fig_05}
\end{figure}

Next, we examine the normalized residual amplitude on several circles in the $(x, y)$-plane. 
Each circle is centered at $\vec{x}_d$, is specified by the radial distance $r$, 
and is parameterized by the azimuthal angle $\phi$ about the center $\vec{x}_d$. 
Since the system (including the source, the given target region of masking, and 
the two masking point-forces) is axial-symmetric about the $x$-axis, we only need to study 
the upper half of each circle.  
\[ \vec{x}(\phi) = \vec{x}_d + r(\cos\phi, \sin\phi, 0), \quad 0 \le \phi \le \pi \]
We consider five semi-circles with radius $r = 1,\, 2,\, 4,\, 8,\, 15$.
Figure \ref{fig_05} plots the normalized residual amplitude 
$A^\text{(NR)}(\vec{x})$ versus azimuthal angle $\phi$ along five semi-circles. 
We compare the results of two masking configurations: 
source + force 1 (left panel); source + the optimal weighted average of forces 1 and 2 
(right panel).
In both panels, for small $r$, the maximum residual amplitude occurs at $\phi = 0$ 
and $\phi = \pi$ (i.e., along the $\vec{x}_d$-direction). 
For moderate values of $r$, the location of maximum residual amplitude shifts to 
$\phi = \pi/2$ (i.e., along the direction perpendicular to $\vec{x}_d$). 
Once again, Figure \ref{fig_05} demonstrates that the optimal weighted average 
of the two point-forces is significantly more effective than either individual point-force
in masking the Gaussian source in the given target region around $\vec{x}_d$. 

\begin{figure}[!h]
\vskip 0.2cm
\begin{center}
\psfig{figure=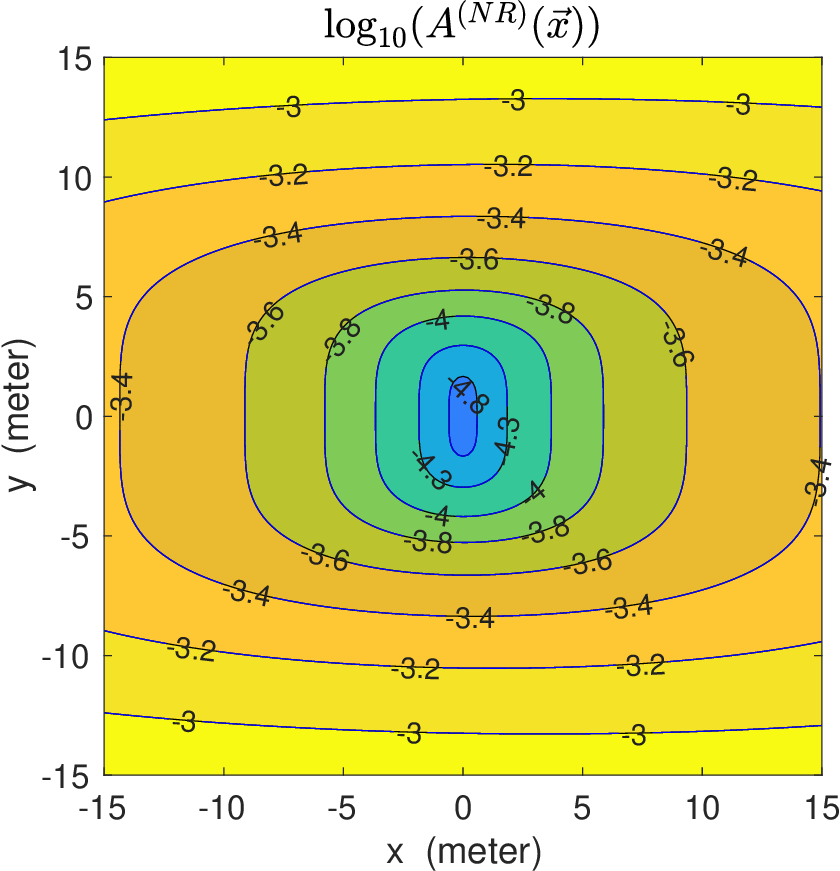, width=2.8in}\qquad
\psfig{figure=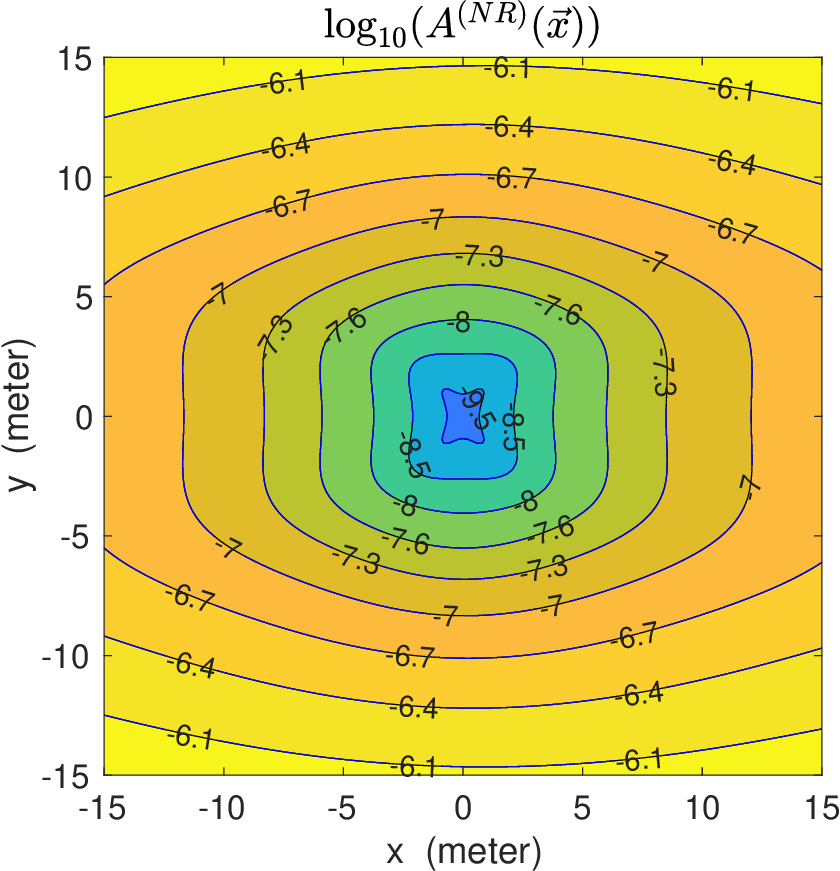, width=2.8in} 
\end{center}
\vskip -0.75cm
\caption{Contours of $\log_{10}\big(A^\text{(NR)}(\vec{x})\big)$
in local coordinates $(x, y)$ about $\vec{x}_d$ where $A^\text{(NR)}(\vec{x})$ 
is the normalized residual amplitude defined in \eqref{resi_amp_def}. 
Left: source + force 1. Right: source + the optimal weighted average of forces 1 and 2.}
\label{fig_06}
\end{figure}

The normalized residual amplitude $A^\text{(NR)}(\vec{x})$ is a function 
$\vec{x} \in \mathbb{R}^3$. Because of the axial-symmetric about the $x$-axis, 
we only need to study $A^\text{(NR)}(\vec{x})$ in the $(x, y)$-plane. 
We write $\vec{x} = \vec{x}_d + (x, y, 0)$ and view $A^\text{(NR)}(\vec{x})$ 
as a function of $(x, y)$, the local coordinates about $\vec{x}_d$ in the $(x, y)$-plane. 
Figure \ref{fig_06} presents contour plots of $\log_{10}\big(A^\text{(NR)}(\vec{x})\big)$ 
in the local coordinates $(x, y)$ about $\vec{x}_d$. 
We compare the results of two masking configurations: 
source + force 1 (left panel); source + the optimal weighted average of forces 1 and 2 
(right panel).

By the design of the masking force, the residual amplitude vanishes at $(x, y) = (0, 0)$ in local coordinates, 
corresponding to $\vec{x}=\vec{x}_d$, the center of the target region in masking. 
At a location displaced radially from $\vec{x}_d$, 
the residual amplitude increases with the radial distance. This increase is notably more 
pronounced along the direction perpendicular to $\vec{x}_d$, which corresponds to the 
vertical direction in Figure \ref{fig_06}.
%

\subsection{Optimal Solution for Masking in a Moderately Large Neighborhood}
The solution given in \eqref{F_wt_avg_2}, together with \eqref{optimal_pt_frc_1}
and \eqref{optimal_pt_frc_2}, describes the optimal set of two point-forces in masking 
the effect of the source in a small neighborhood around $\vec{x}_d$, 
the center of the region where possible detecting sensors are located. 
Mathematically, this solution is the optimal set of two point-forces in masking 
in an infinitesimal neighborhood. 
The optimization analysis is based on the Taylor approximation around $\vec{x}_d$.
When masking in a moderately large neighborhood of $\vec{x}_d$, however, 
this solution is no longer truly optimal. 
It is possible to adjust the two-force configuration to further reduce 
the maximum residual amplitude over the finite target region. 

Recall that in the small-neighborhood optimal solution given in \eqref{optimal_pt_frc_1}
with \eqref{optimal_pt_frc_2} and \eqref{F_wt_avg_2}, the parameters are derived 
based on the Taylor approximation around $\vec{x}_d$. The specific strategy is that 
we adjust the parameters to make 
\begin{equation}
\begin{dcases}
u(\vec{x}, t) \Big|_{\vec{x} = \vec{x}_d} = 0  \\[1ex]
\nabla u(\vec{x}, t) \Big|_{\vec{x} = \vec{x}_d} = 0
\end{dcases} 
\label{set_u_Du_0}
\end{equation}
The two conditions in \eqref{set_u_Du_0} are equivalent to minimizing the residual amplitude
over an infinitesimal neighborhood around $\vec{x}_d$. 
Among the parameters in the two-force configuration, $\varphi_m^{(1)}$ and 
$\varphi_m^{(2)}$ are the phase shifts, respectively, in the two oscillating point-forces. 
When masking in a small neighborhood, $\varphi_m^{(1)}$ and $\varphi_m^{(2)}$ are 
determined from the condition $u(\vec{x}, t) \big|_{\vec{x} = \vec{x}_d} = 0$ and are 
given in \eqref{optimal_pt_frc_1} and \eqref{optimal_pt_frc_2}. 
When masking in a moderately large neighborhood, minimizing the residual amplitude 
over the finite region is in general inconsistent with the constraint 
$u(\vec{x}, t) \big|_{\vec{x} = \vec{x}_d} = 0$. 
Thus, in the new optimization for the finite region, we set parameters $\varphi_m^{(1)}$ 
and $\varphi_m^{(2)}$ free while keeping other parameters given in 
\eqref{optimal_pt_frc_1} and \eqref{optimal_pt_frc_2}. 

Let $A^\text{(NR)}(\vec{x}; \varphi_m^{(1)}, \varphi_m^{(2)} )$ 
be the normalized residual amplitude after masking defined in \eqref{resi_amp_def}. 
In the extended notation of $A^\text{(NR)}$, we explicitly include its dependence on 
$(\varphi_m^{(1)}, \varphi_m^{(2)})$, the two free parameters in our the new optimization. 
Let $B(\vec{x}_d, r_d)$ be the target region for masking. The size of the region 
is described by the radius $r_d$. We consider the spatial maximum of 
$A^\text{(NR)}(\vec{x}; \varphi_m^{(1)}, \varphi_m^{(2)} )$ 
over the region $B(\vec{x}_d, r_d)$. 
\begin{equation}
E_\text{max}(r_d; \varphi_m^{(1)}, \varphi_m^{(2)} ) \equiv 
\max_{\vec{x} \in B(\vec{x}_d, r_d)} 
A^\text{(NR)}(\vec{x}; \varphi_m^{(1)}, \varphi_m^{(2)}) 
\label{E_max_def} 
\end{equation}
In the new optimization, we minimize $E_\text{max}(r_d; \varphi_m^{(1)}, \varphi_m^{(2)}) $ 
over $(\varphi_m^{(1)}, \varphi_m^{(2)})$. 
\begin{align}
(\varphi_m^{(1)}, \varphi_m^{(2)})^\text{(opt)}(r_d) \equiv 
 \argmin_{(\varphi_m^{(1)}, \varphi_m^{(2)})} 
 E_\text{max}(r_d; \varphi_m^{(1)}, \varphi_m^{(2)} )
 \label{opt_at_rd} 
\end{align}
For a finite region, we expect that the optimal solution varies with the radius of the region. 
Optimization problem \eqref{opt_at_rd} is solved numerically for each value of $r_d$. 
Mathematically, the expressions of $(\varphi_m^{(1)}, \varphi_m^{(2)})$ given in 
\eqref{optimal_pt_frc_1} and \eqref{optimal_pt_frc_2} represent 
the true optimal solution in the limit of an infinitesimally small radius $r_d$. That is,
\[ 
\Big[(\varphi_m^{(1)}, \varphi_m^{(2)}) \text{ given in } 
\text{\eqref{optimal_pt_frc_1} and \eqref{optimal_pt_frc_2}}\Big] = 
\lim_{r_d \rightarrow 0_+} (\varphi_m^{(1)}, \varphi_m^{(2)})^\text{(opt)}(r_d) 
\]
For a small $r_d$, the optimal solution is approximated by its limit as $r_d \rightarrow 0$. 
\[ (\varphi_m^{(1)}, \varphi_m^{(2)})^\text{(opt)}(r_d ) 
\approx (\varphi_m^{(1)}, \varphi_m^{(2)})^\text{(opt)}(0_+) \quad 
\text{for small $r_d$} \]
We compare two slightly different designs of two-force configuration for masking 
in a given finite region $B(\vec{x}_d, r_d)$: 
the first one is the optimal solution for an infinitesimal region; 
the second one is the result of the new optimization specifically tailored to the given finite region. 
The only difference between the two designs is in parameters 
$(\varphi_m^{(1)}, \varphi_m^{(2)})$. 
\begin{itemize}
\item Design 1: we use $(\varphi_m^{(1)}, \varphi_m^{(2)})^\text{(opt)}(0_+) $. \\
The maximum residual amplitude over the given target region $B(\vec{x}_d, r_d)$ is 
\begin{align}
 & E_\text{max}^\text{(opt, $0_+$)}(r_d)
 \equiv E_\text{max}(r_d; (\varphi_m^{(1)}, \varphi_m^{(2)})^\text{(opt)}(0_+) ) 
 \nonumber \\
& \qquad = \max_{\vec{x} \in B(\vec{x}_d, r_d)} 
A^\text{(NR)}(\vec{x}; (\varphi_m^{(1)}, \varphi_m^{(2)})^\text{(opt)}(0_+) ) 
\label{E_max_1}
\end{align}
\item Design 2: we use $(\varphi_m^{(1)}, \varphi_m^{(2)})^\text{(opt)}(r_d)$. \\
The maximum residual amplitude over the given target region $B(\vec{x}_d, r_d)$ is 
\begin{align}
 & E_\text{max}^\text{(opt, $r_d$)}(r_d)
 \equiv E_\text{max}(r_d; (\varphi_m^{(1)}, \varphi_m^{(2)})^\text{(opt)}(r_d) ) 
 \nonumber \\
& \qquad = \min_{(\varphi_m^{(1)}, \varphi_m^{(2)})} \max_{\vec{x} \in B(\vec{x}_d, r_d)} 
A^\text{(NR)}(\vec{x}; \varphi_m^{(1)}, \varphi_m^{(2)} ) 
\label{E_max_2}
\end{align}
\end{itemize}
It follows from the additional minimization in \eqref{E_max_2} that  
\[ E_\text{max}^\text{(opt, $r_d$)}(r_d) \le 
E_\text{max}^\text{(opt, $0_+$)}(r_d) \]
\begin{figure}[!h]
\vskip 0.25cm
\begin{center}
\psfig{figure=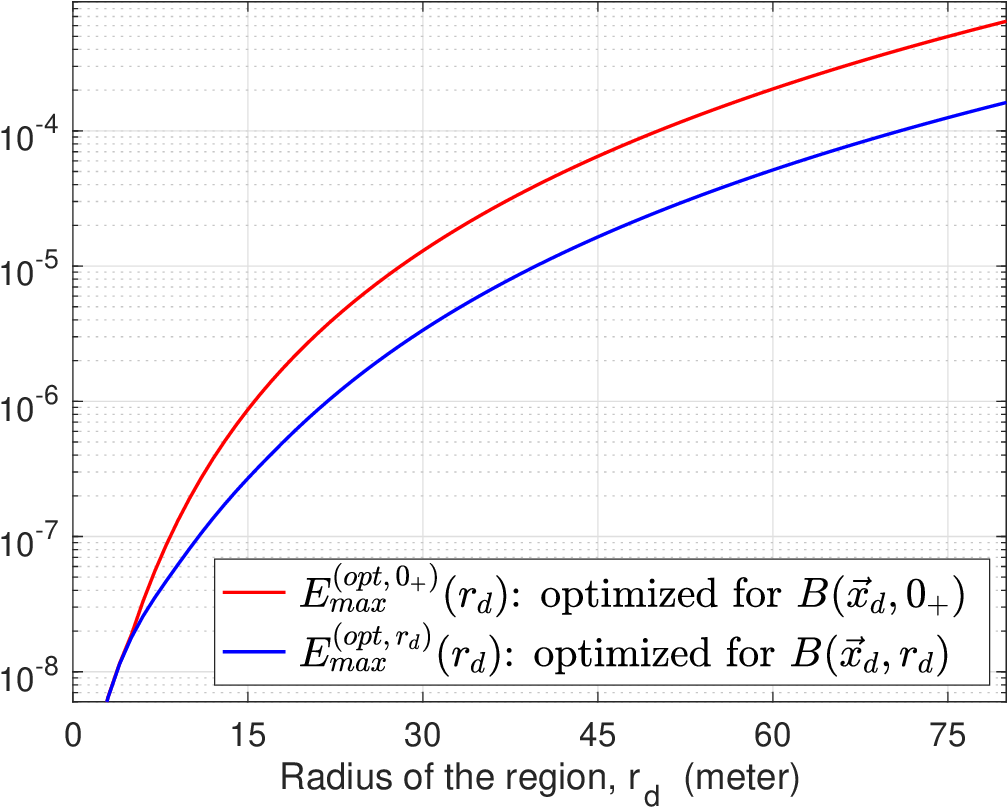, width=2.85in}\qquad
\psfig{figure=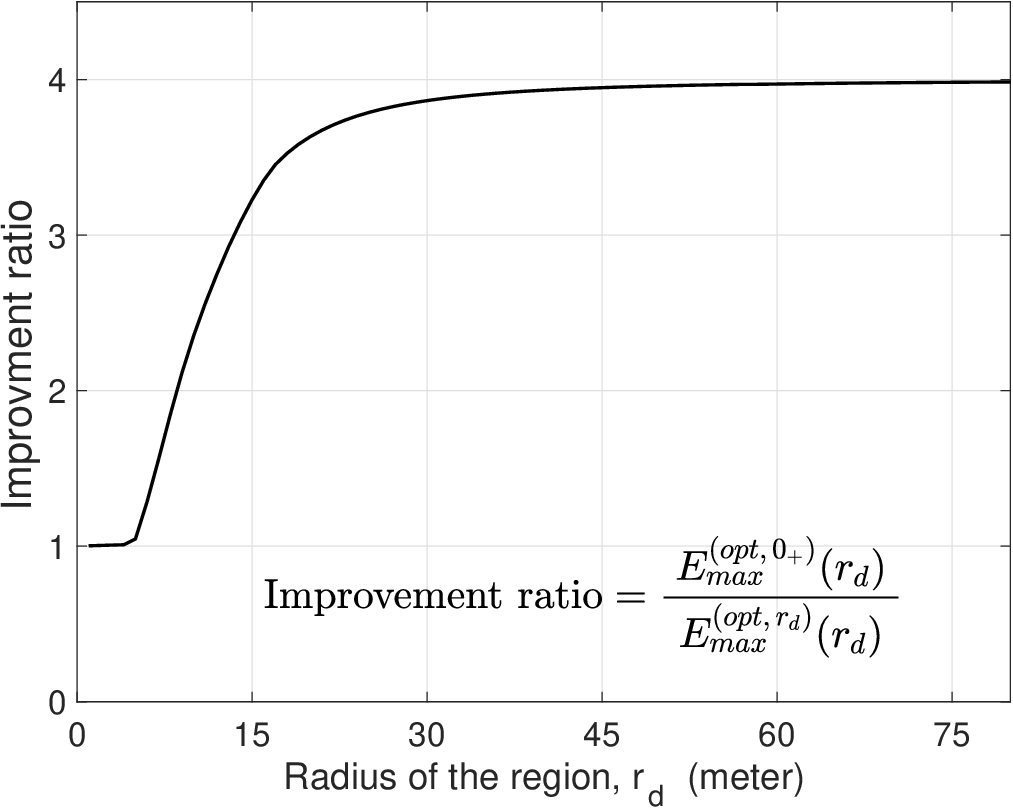, width=2.85in} 
\end{center}
\vskip -0.75cm
\caption{Maximum residual amplitudes of the two designs defined in \eqref{E_max_1}
and \eqref{E_max_2}. }
\label{fig_07}
\end{figure}
The left panel of Figure \ref{fig_07} compares $E_\text{max}^\text{(opt, $0_+$)}(r_d)$
and $E_\text{max}^\text{(opt, $r_d$)}(r_d)$ as functions of $r_d$, 
the radius of the target region. 
For a small $r_d$ (i.e., masking in a small neighborhood of $\vec{x}_d$), 
the two designs have virtually the same maximum residual amplitudes. 
For a moderately large $r_d$ (i.e., masking in a moderately large neighborhood 
of $\vec{x}_d$), the numerically optimized 
$(\varphi_m^{(1)}, \varphi_m^{(2)})^\text{(opt)}(r_d)$ outperforms 
the small-neighbor approximation $(\varphi_m^{(1)}, \varphi_m^{(2)})^\text{(opt)}(0_+) $. 

The right panel of Figure \ref{fig_07} plots the ratio of $E_\text{max}^\text{(opt, $0_+$)}(r_d)$
to $E_\text{max}^\text{(opt, $r_d$)}(r_d)$. 
For $r_d > 15\,\text{m}$, the maximum residual amplitude after masking is improved 
by a factor of 3 or more. The improvement shown is the additional gain achieved by switching from 
the small-neighbor approximation (design 1) to the optimal solution that is specifically tailored 
to the given target region $B(\vec{x}_d, r_d)$ (design 2). 

\begin{figure}[!h]
\vskip 0.25cm
\begin{center}
\psfig{figure=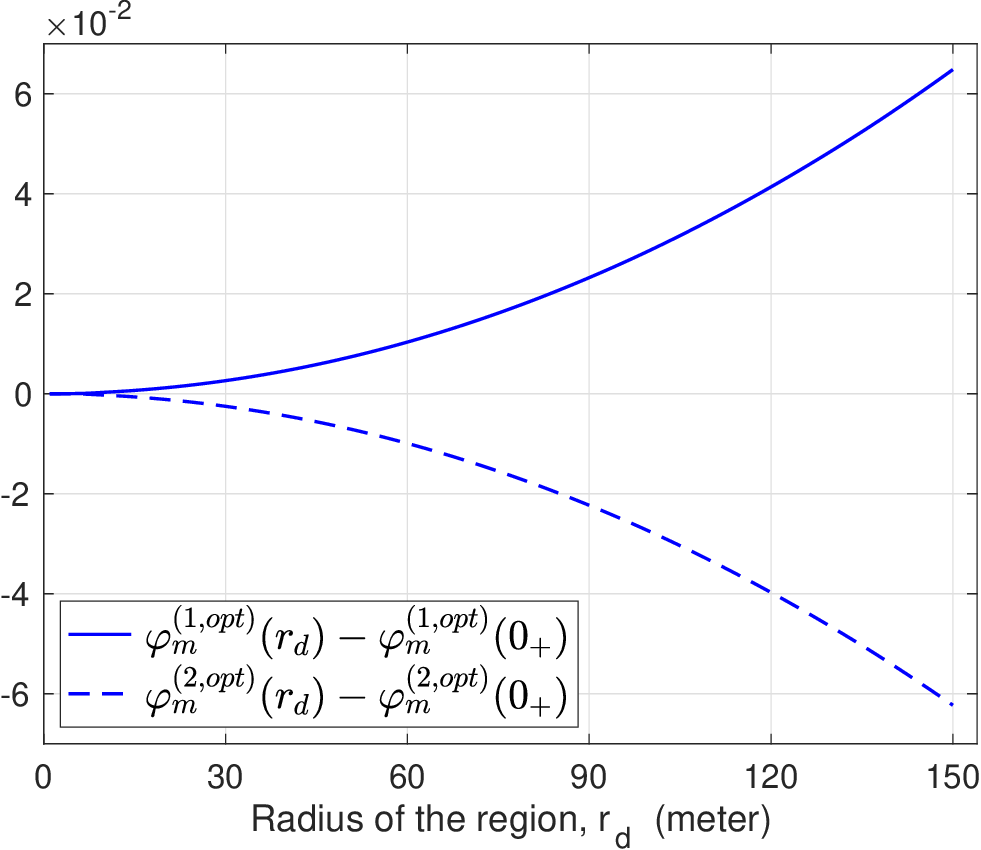, width=3in}
\end{center}
\vskip -0.75cm
\caption{$\varphi_m^\text{(1, opt)}(r_d)$ and $\varphi_m^\text{(2, opt)}(r_d)$
defined in \eqref{opt_at_rd} as functions of $r_d$. }
\label{fig_08}
\end{figure}

We investigate how the optimal phase shifts $(\varphi_m^{(1)}, \varphi_m^{(2)})^\text{(opt)}(r_d )$
vary with the radius $r_d$ of the target region. 
To highlight the deviations from the small-neighborhood approximations, in Figure \ref{fig_08}, 
we plot $ \big(\varphi_m^\text{(1, opt)}(r_d )- \varphi_m^\text{(1, opt)}(0_+)\big)$ 
and $ \big(\varphi_m^\text{(2, opt)}(r_d )- \varphi_m^\text{(2, opt)}(0_+)\big)$ vs $r_d$. 
The vertical axis of Figure \ref{fig_08} is in units of radians. In an oscillating force,
a phase shift of  $6.28\times 10^{-2}$ radian corresponds to $1/100$ of an oscillation cycle. 

\begin{figure}[!h]
\vskip 0.2cm
\begin{center}
\psfig{figure=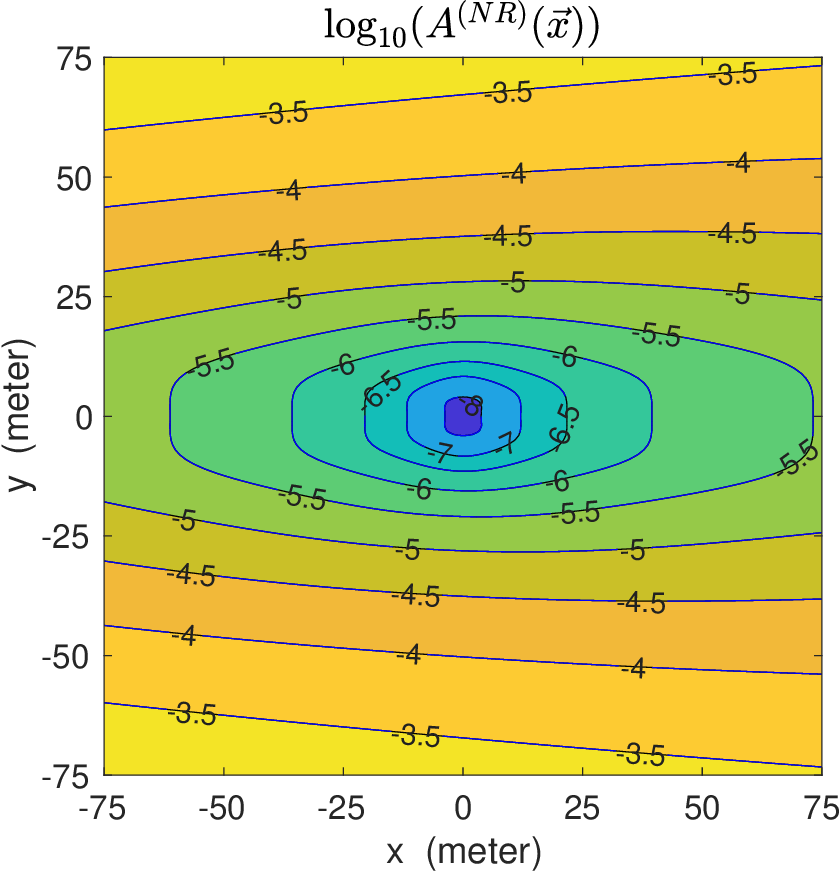, width=2.8in}\qquad
\psfig{figure=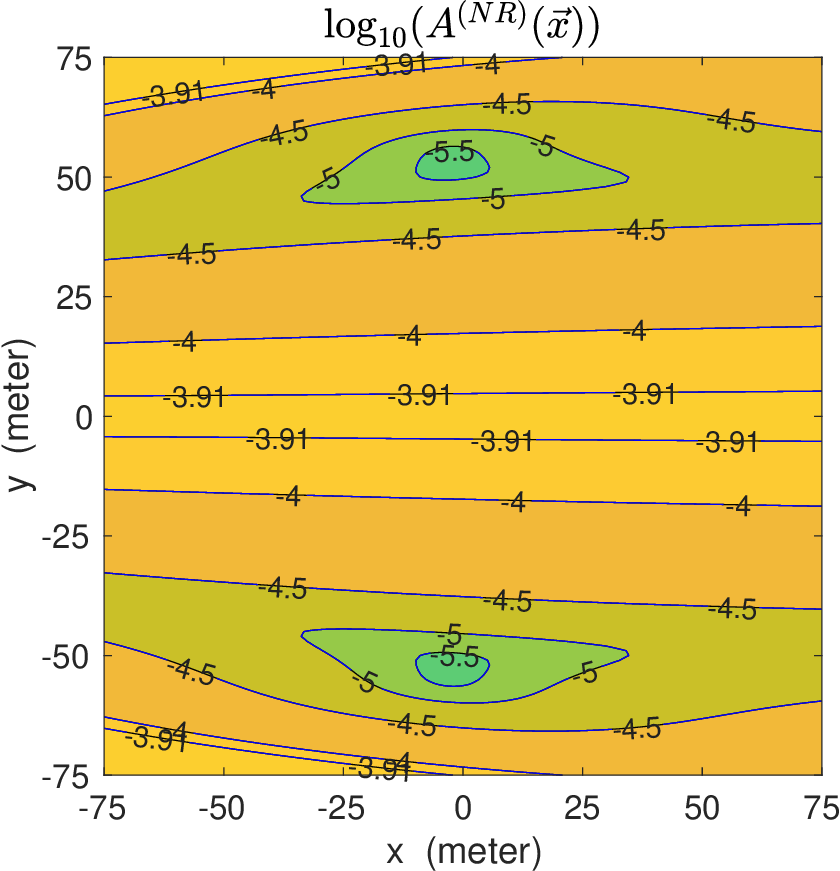, width=2.8in} 
\end{center}
\vskip -0.75cm
\caption{Contours of $\log_{10}\big(A^\text{(NR)}(\vec{x})\big)$ in a moderately 
large region $B(\vec{x}_d, 75\,\text{m})$. 
Left: performance of the two-force configuration optimized for $B(\vec{x}_d, 0_+)$.
Right: performance of the two-force configuration optimized for $B(\vec{x}_d, 75\,\text{m})$.}
\label{fig_09}
\end{figure}

We compare the performance of the two designs for masking the effect of a Gaussian source
in a target region of radius $r_d = 75\, \text{meters}$. Recall that in simulations, 
the center of the target region is $| \vec{x}_d | = 750\, \text{meters}$ 
away from the acoustic source. Relative to this distance, a ball of radius $r_d = 75\, \text{meters}$ 
is a fairly large region. 
Figure \ref{fig_09} plots the contours of $\log_{10}\big(A^\text{(NR)}(\vec{x})\big)$
in the local coordinates $(x, y)$ about $\vec{x}_d$ where $A^\text{(NR)}(\vec{x})$ 
is the normalized residual amplitude after applying a two-force masking configuration, 
defined in \eqref{resi_amp_def}. 
In the left panel, the two-force configuration is optimized for an infinitesimal region but 
is applied to the finite region of radius $r_d = 75\, \text{meters}$. 
Since the optimization is based on the Taylor approximation about $(0, 0)$, 
it is guaranteed that the residual amplitude vanishes at $(0, 0)$. 
In the right panel the two-force configuration is optimized numerically specifically 
for the given finite region. In the new optimization, we no longer have zero 
residual amplitude at $(0, 0)$. In return, we are able to achieve 
a lower maximum residual amplitude over the finite region. 

The numerical optimization for the given finite region is with respect to 
$(\varphi_m^{(1)}, \varphi_m^{(2)})$ while all other parameters are fixed 
 at the values of small-neighborhood approximation given in 
 \eqref{optimal_pt_frc_1} and \eqref{optimal_pt_frc_2}.
Before concluding this section, we point out that the maximum residual amplitude
can be further reduced by expanding the set optimization variables. 
For example, if we allow both the phase shifts $(\varphi_m^{(1)}, \varphi_m^{(2)})$ 
and the amplitudes $(a_m^{(1)}, a_m^{(2)})$ to vary freely in optimization, 
we can obtain a slightly lower maximum residual amplitude. 
This underscores the inherent complexity of designing optimal inference forces 
for masking an acoustic source in a specified target region of possible detection sensors. 
Even in the relatively simple setting of two-force configurations, a full optimization 
without any theoretical constraints involves 10 variables in total: 
3 spatial coordinates for each force location, and 1 amplitude and 1 phase shift 
for each point-force. 
%

\section{A Numerical Method for Problems Lacking Spherical Symmetry} \label{num_method}
We present an efficient numerical method for solving the 3D forced wave equation 
in the unbounded domain $\mathbb{R}^3$. The objective is to accommodate the situation 
where the forcing term $F(\vec{x}, t) = p(\vec{x}) \sin(\omega t) $ lacks 
spherical symmetry in space. 
The spherically symmetric case was analytically resolved in 
Section \ref{wav_sph_sym}.
The numerical method is based on the concept of Green's functions and leverages
the analytical solutions obtained in Section \ref{wav_sph_sym}.

Let $h(s)$ be a smooth function of unit spatial scale satisfying 
i) $ Z_h \equiv \int_0^{+\infty} 4\pi s^2 h(s) ds > 0$ and 
ii) $h(s)$ decaying to 0 quickly as $s \rightarrow +\infty$. 
For example, $h(s)$ may be the standard normal density $\rho_{N(0, 1)}(s)$ or 
the truncated sinc function $\text{sinc}^{(tc)}(s)$ introduced in Section \ref{wav_sph_sym}. 
We use $h(s)$ to build a smoothed version of the Dirac delta function 
$\delta(\vec{x}-\vec{x}_0)$. Let 
\begin{equation}
\psi(\vec{x}; \vec{x}_0, d) \equiv \frac{1}{Z_h d^3 } h(\frac{1}{d}|\vec{x}-\vec{x}_0|) 
\label{psi_def}
\end{equation}
$\psi(\vec{x}; \vec{x}_0, d)$ has the properties below.
\begin{itemize}
\item $\psi(\vec{x}; \vec{x}_0, d)$ has spherical symmetry about $\vec{x}_0$: 
\[ \psi(\vec{x}; \vec{x}_0, d) = \psi(r; 0, d) = \frac{1}{Z_h d^3 } h(\frac{r}{d}), 
\quad r \equiv | \vec{x}-\vec{x}_0 | \] 
\item $\displaystyle \int_{\mathbb{R}^3} \psi(\vec{x}; \vec{x}_0, d) d \vec{x} 
= \frac{1}{Z_h} \int_0^{\infty} 4\pi s^2 h(s) ds = 1$. 
\item The effective support of $\psi(\vec{x}; \vec{x}_0, d) $ is proportional to $d$. 
\item $\displaystyle \lim_{d \rightarrow 0_+} \psi(\vec{x}; \vec{x}_0, d) 
= \delta(\vec{x}-\vec{x}_0)$. 
\end{itemize}
Consider a general sinusoidally oscillating force $F(\vec{x}, t) = p(\vec{x}) \sin(\omega t)$.
We solve the general case where the spatial profile $p(\vec{x})$ lacks spherical symmetry. 
Mathematically, we can always write $p(\vec{x})$ as an integral 
superposition of the delta function. 
\[ p(\vec{x}) = \int_{\mathbb{R}^3} p(\vec{x}') \delta(\vec{x}-\vec{x}') d \vec{x}' \]
For small $d$, function $\psi(\vec{x}; \vec{x}_j, d)$ can be viewed as a smoothed version 
of $\delta(\vec{x}-\vec{x}_j)$. 
Thus, we approximate $p(\vec{x})$ using a summation superposition of 
$\psi(\vec{x}; \vec{x}_j, d)$. 
\begin{equation}
p(\vec{x}) \approx \sum_{j=1}^N a_j \psi(\vec{x}; \vec{x}_j, d)  
\label{px_approx}
\end{equation}
where $\{\vec{x}_j \}$ is a set of discrete points representing the effective support 
of $p(\vec{x})$ and $\{a_j \}$ is the optimal set of coefficients that minimizes the error of 
approximation in \eqref{px_approx}. 
Since the 3D forced wave equation is linear, we only need to solve the case of 
$F(\vec{x}, t) = \psi(\vec{x}; \vec{x}_j, d) \sin(\omega t)$, which is spherically symmetric 
after shifting by $\vec{x}_j$. 
\[ \psi(\vec{x}; \vec{x}_j, d) \sin(\omega t) = \frac{1}{Z_h }\cdot \frac{1}{d^3 } h(\frac{r}{d}) \sin(\omega t), 
\quad r \equiv | \vec{x}-\vec{x}_j | \] 
Let $u\big(r, t; \{h(s\}, d, \omega \big)$
be the quasi-steady periodic solution for $F(r, t) = \frac{1}{d^3 } h(\frac{r}{d}) \sin(\omega t)$, 
whose close-form analytical expression was derived in Section \ref{wav_sph_sym}. 
Based on the principle of superposition, the solution for the general 
$F(\vec{x}, t) = p(\vec{x}) \sin(\omega t)$ is approximately 
\begin{equation}
\boxed{\quad
u(\vec{x}, t) \approx \sum_{j=1}^N \frac{a_j}{Z_h} u\big(|\vec{x}-\vec{x}_j|, t; \{h(s)\}, d, \omega \big) 
\quad}
\label{num_sol_u}
\end{equation}
In the case where $h(s) = \text{sinc}^{(tc)}(s)$, the solution 
$u\big(r, t; \{h(s)\}, d, \omega \big)$ is given in \eqref{u_sinc_QSS}; 
in the case where $h(s) = \rho_{N(0, 1)}(s)$, the solution 
$u\big(r, t; \{h(s)\}, d, \omega \big)$ is given in \eqref{u_normal_QSS}. 
%

\section{Conclusions}
This paper presents a theoretical and computational framework for designing 
interference signals that minimize the detectability of a given acoustic source 
in a known target region where possible sensors are located. 
We derived analytical quasi-steady periodic solutions to the three-dimensional (3D) forced 
wave equation under spherical symmetry in several canonical cases. 
We then extended the analysis to practical masking scenarios, including 
(1) self-masking where an acoustic source, with certain spatial forcing profile and 
certain spatial scale relative to the wavelength, cancel itself outside its effective forcing core; 
and (2) localized masking in a known target region by placing one or 
two interference forces near the acoustic source. In general, the given acoustic source
is at a moderate distance away from the target region, not right next to it. 
In the case of one-force configuration or two-force configuration, respectively, 
the optimal solution for masking in a small target region is analytically derived based on 
the Taylor expansion. For two-force masking, the optimal solution 
makes both the residual amplitude and the gradient of residual amplitude vanish 
at the center of the target region, which mathematically minimizes the maximum 
residual amplitude in a small region. The masking performances of the one-force configuration
and the two-force configuration are tested and compared in numerical simulations. 
When masking in a moderately large region, the analytical solution based on 
the Taylor expansion is no longer the optimal. 
In this case, we carried out numerical optimization specifically tailored to the finite target region. 
The maximum residual amplitude is further reduced by a factor of 3 or more 
in the numerical optimization. 

In the case of a more general forcing that lacks spherical symmetry, we presented 
a numerical method for solving the 3D forced wave equation in the 
unbounded 3D domain. 
The numerical method is based on representing the general force approximately 
as a linear combination of spherically symmetric kernels, distributed over the effective support 
of the spatial forcing profile, each kernel being a smoothed version of the Dirac delta function. 
The solution excited by each kernel 
has spherical symmetry about its center and is analytically derived. The solution excited by 
the general force is approximated by a superposition of these kernel solutions. This numerical 
approach is completely adaptive and very efficient; 
there is no 3D numerical grid involved in computation. 

Potential applications of this work include secure communication, submarine stealth, 
and defense against acoustic surveillance. Future work will explore further optimization strategies, 
extend the method to incorporate more general forcing profiles, and apply these techniques
to real-world problems in wave control and signal masking.

\vskip 2cm 
\noindent{\bf Acknowledgments}\\
HZ acknowledges the support provided by the Undersea Warfare Academic Group at the Naval Postgraduate School.

\end{document}